\documentclass{amsart}

\newtheorem{theorem}{Theorem}[section]
\newtheorem{lemma}[theorem]{Lemma}
\newtheorem{proposition}[theorem]{Proposition}
\newtheorem{corollary}[theorem]{Corollary}

\theoremstyle{definition}
\newtheorem{definition}[theorem]{Definition}

\theoremstyle{remark}
\newtheorem{remark}[theorem]{Remark}

\numberwithin{equation}{section}

\newcommand{\abs}[1]{\lvert#1\rvert}

\newcommand{\M}{\mathcal{M}}

\newcommand{\p}{\mathbb{P}}
\newcommand{\I}{\mathbb{I}}

\begin{document}
\title[Pure state transformations]{Pure state
transformations induced by linear operators}
\author{L E Labuschagne}
\address{Department of Mathematical Sciences, P.O.Box 392,
University of South Africa, 0003 Pretoria, South Africa}
\email{labusle@unisa.ac.za}


\thanks{Part of this research was conducted with the support of a joint
KBN-NRF grant under the Poland-South Africa cooperation agreement.
The author would also like to acknowledge some inspiring
discussions with prof WA Majewski.}

\begin{abstract}
We generalise Wigner's theorem to its most general form possible for
$B(\mathfrak{h})$ in the sense of completely characterising those
vector state transformations of $B(\mathfrak{h})$ that appear as restrictions
of duals of linear operators on $B(\mathfrak{h})$. We then use this result to
similarly characterise all pure state transformations of general
$C^\ast$-algebras that appear as restrictions of duals of linear
operators on the underlying algebras. This result may variously be
interpreted as either a non-commutative Banach-Stone theorem, or (in the
bijective case) a pure state based description of Wigner symmetries.
These results extend the work of Shultz \cite{Sh} (who considered
only the case of bijections), and also complements and completes
the investigation of linear maps with pure state preserving
adjoints begun in \cite{LMa}.
\end{abstract}

\maketitle

\section{Preliminaries}

Notation in this paper will be based on that of \cite{BRo} and
\cite{KRi}. We briefly review the essentials. In general
$\mathcal{A, B}$ will denote $C^\ast$-algebras with
$\mathcal{A}_{sa}$ denoting the self-adjoint portion of a given
$C^\ast$-algebra $\mathcal{A}$. Unless otherwise stated we will
generally assume our algebras to be unital with $\I$ being used to
denote the unit. In this context a linear map $\varphi:
\mathcal{A} \rightarrow \mathcal{B}$ is called a Jordan
$\ast$-homomorphism if it preserves both the Jordan product and
adjoints, ie. if $\varphi(AB + BA) = \varphi(A)\varphi(B) +
\varphi(B)\varphi(A)$ and $\varphi(A^\ast) =  \varphi(A)^\ast$ for
all $A, B \in \mathcal{A}$.

Given a $C^\ast$-algebra $\mathcal{A}$, the state space will be denoted by
$\mathcal{S}_{\mathcal A}$ with $\omega, \rho$ denoting typical states on
$\mathcal{A}$. At a slight variance with the more common practice we will use
$\mathcal{P}_{\mathcal A}$ to denote merely the set of pure states of
$\mathcal{A}$, rather than the $\mathrm{weak}^\ast$-closure of this set. If
$B(\mathfrak{h})$ is the space of bounded operators on some Hilbert space
$\mathfrak{h}$, a representation $\pi: \mathcal{A} \rightarrow
B(\mathfrak{h})$ of $\mathcal{A}$ as bounded operators on $\mathfrak{h}$ will
be denoted by $(\pi, \mathfrak{h})$. If indeed this is the representation
canonically engendered by some state $\omega$ via the Gelfand-Neumark-Segal
construction, this will be indicated by a suitable subscript. If now
$\mathcal{A}$ is a concrete $C^\ast$-algebra on some Hilbert space
$\mathfrak{h}$, a vector state $\mathcal{A} \rightarrow \mathbb{C} :
A \mapsto (Ax, x)$ yielded by some norm-one vector $x \in \mathfrak{h}$ will
be denoted by $\omega_x$.

When focussing specifically on von Neumann algebras, we will employ the
notation $\M_1, \M_2, ...$ for these algebras. Given a von Neumann algebra
$\M$, its projection lattice will be denoted by $\mathbb{P}_\M$. In this
context an \emph{orthoisomorphism} between two von Neumann algebras is
understood to be a bijection between their respective projection lattices
which preserves orthogonality of projections. More specifically given two von
Neumann algebras $\M_1$ and $\M_2$, $\psi: \mathbb{P}_{\M_1} \rightarrow
\mathbb{P}_{\M_2}$ is an orthoisomorphism if it injectively maps
$\mathbb{P}_{\M_1}$ onto $\mathbb{P}_{\M_2}$ and also satisfies the condition
that $$EF = 0 \quad \mbox{if and only if} \quad \psi(E)\psi(F) = 0$$ for all
$E, F \in \mathbb{P}_{\M_1}$.

Given two $C^\ast$-algebras $\mathcal{A}$ and $\mathcal{B}$ each corresponding
to some quantum mechanical system, it is natural to consider two such algebras
to be \emph{physically equivalent} if in some phenomenological sense they are
the same. Based on the fact that the essential physical information of a
quantum mechanical system is encoded in its set of possible states, Emch
defined two representations $(\pi_1, \mathfrak{h}_1)$ and $(\pi_2,
\mathfrak{h}_2)$ of a $C^\ast$-algebra $\mathcal{A}$ to be physically
equivalent if indeed
$$\{\omega \circ \pi_1 | \omega \in \mathcal{S}_{\pi_1({\mathcal A})}\} =
\{\omega \circ \pi_2 | \omega \in \mathcal{S}_{\pi_2({\mathcal A})}\}$$
(see \cite{E}, p 107). Taking this idea a step further, it seems natural to
consider $\mathcal{A}$ and $\mathcal{B}$ to be physically equivalent if their
respective state spaces are structurally identical. More precisely
$\mathcal{A}$ and $\mathcal{B}$ are called physically equivalent if we can
find an affine $\mathrm{weak}^\ast$-continuous bijection $\varphi^\sharp$ from
the state space of $\mathcal{B}$ onto that of $\mathcal{A}$. Such a map
$\varphi^\sharp$ will then be considered to be a map which \emph{preserves all
relevant physical information} in a change of contexts from $\mathcal{A}$ to
$\mathcal{B}$. At the level of von Neumann algebras physical equivalence of
two von Neumann algebras $\M_1$ and $\M_2$ may similarly be defined in terms
of the structural similarity of their respective normal state spaces. In this
setting the mappings which \emph{preserve all relevant physical information}
in a change of contexts from $\M_1$ to $\M_2$ (thereby ensuring the physical
equivalence of these algebras) may be identified with affine bijections from
the normal state space of $\M_2$ onto that of $\M_1$.

Our first task will be to show that various alternative ways of
formalising this notion are in fact equivalent and invariably tend
to lead one to a Jordan $\ast$-isomorphism between the algebras.
To afford such results we need the following lemmas. The second
lemma is of course well-known for $\ast$-isomorphisms but for our
purposes the Jordan case needs to be included. Since both these
lemmas will be needed again later on, we show how Lemma 1.1 may be
used to reduce Lemma 1.2 to the $\ast$-isomorphism case.

\begin{lemma} Let $\M_1$ and $\M_2$ be von Neumann algebras and let
$\varphi$ be a Jordan $\ast$-isomorphism from $\M_1$ onto $\M_2$. Then
there exists a central projection $E \in \M_1 \cap \M_1'$ such that
$\varphi(E) = F$ is a central projection in $\M_2 \cap \M_2'$ for which
$F\varphi(\cdot)F = \varphi_F$ defines a $\ast$-isomorphism from $(\M_1)_E$
onto $(\M_2)_F$ and $(\I-F)\varphi(\cdot)(\I-F) = \varphi_{(\I-F)}$ a
$\ast$-antiisomorphism from $(\M_1)_{(\I-E)}$ onto $(\M_2)_{(\I-F)}$.
\end{lemma}

\begin{proof} The proof follows from an application of Kadison's result
\cite[3.2.2]{BRo} regarding the existence of a projection $F \in
\M_2 \cap (\M_2)'$ such that $\varphi_F$ is a $\ast$-homomorphism
and $\varphi_{(\I-F)}$ a $\ast$-antimorphism from $\M_1$ onto
$(\M_2)_F$ and $(\M_2)_{(\I-F)}$ respectively.
\end{proof}

\begin{lemma} Let $\M_1$ and $\M_2$ be von Neumann algebras and let $\varphi$
be Jordan $\ast$-isomorphism from $\M_1$ onto $\M_2$. Then $\varphi$ is
$\sigma$-weakly and $\sigma$-strongly continuous.
\end{lemma}

\begin{proof}
Let $E$ be as in the preceding lemma and suppose that $\{A_{\lambda}\}$
converges $\sigma$-weakly (alternatively $\sigma$-strongly) to $A_0$ in
$\M_1$. But then $\{A_{\lambda}E\}$ and $\{A_{\lambda}(\I-E)\}$ converge
$\sigma$-weakly (alt. $\sigma$-strongly) to $A_0E$ and $A_0(\I-E)$
respectively \cite[2.4.1 \& 2.4.2]{BRo}. By for example \cite[2.4.23]{BRo}
applied to the previous lemma it follows that each of $\varphi_F$
and $\varphi_{(\I-F)}$ is $\sigma$-weakly (alt. $\sigma$-strongly) continuous
on $(\M_1)_E$ and $(\M_1)_{(\I-E)}$ respectively. Consequently
$\{\varphi(A_{\lambda})F\}$ and $\{\varphi(A_{\lambda})(\I-F)\}$ converge
$\sigma$-weakly (alt. $\sigma$-strongly) to $\varphi(A_0)F$ and
$\varphi(A_0)(\I-F)$ respectively. Clearly $\{\varphi(A_{\lambda})\}$ will
then converge $\sigma$-weakly (alt. $\sigma$-strongly) to $\varphi(A_0)$.
\end{proof}

Besides establishing the framework for the results that follow,
the following collection of known results has a lot of
philosophical significance. From the point of view of the
preservation of relevant physical information in a change of
contexts from one algebra to another, this array of results
suggests that in the category of von Neumann algebras not just the
structural similarity of the state spaces (respectively normal
state spaces) will guarantee the physical equivalence of two von
Neumann algebras, but also the similarity of any one of their
metric structures, their order structures, and (provided we
exclude the case of the $2\times2$ matrices) even the similarity
of their quantum propositional calculi. This result may also be
interpreted in terms of Wigner-symmetries. In elementary quantum
mechanics the states of a system may be taken to correspond to the
one-dimensional subspaces of some Hilbert space $\mathfrak{h}$.
Given two unit vectors $x, y \in \mathfrak{h}$, the transition
probability between the corresponding states $\mathrm{span}\{x\}$
and $\mathrm{span}\{y\}$ is defined to be $\abs{(x, y)}^2$. In
this context Wigner \cite{W} defined a symmetry to be a bijection
on these states which preserves the transition probabilities
between states, before going on to show that the class of all such
symmetries corresponds canonically to the class consisting of both
unitarily implemented $\ast$-automorphisms and anti-unitarily
implemented $\ast$-anti-automorphisms on $B(\mathfrak{h})$ (or
equivalently the class of all Jordan $\ast$-automorphisms on
$B(\mathfrak{h})$ - see for example p 210 and Example 3.2.14 of
\cite{BRo}). For a fuller account of the connection between Wigner
symmetries and Jordan $\ast$-automorphisms the reader is referred
to either p 210 of \cite{BRo} or p 150 of \cite{E}. In the light
of the above discussion it therefore makes sense to consider
Jordan $\ast$-isomorphisms from one $C^\ast$-algebra onto another
as some sort of Wigner symmetry. With this in mind, the following
result may be interpreted as saying that not just its action on
the relevant state spaces serves to identify a transformation as a
Wigner symmetry, but its action on various other structures as
well.

\begin{theorem} Let $\M_1$ and $\M_2$ be von Neumann algebras. The set of
all Jordan $\ast$-isomorphisms from $\M_1$ onto $\M_2$ is in a one
to one correspondence with each of the following in the sense of
each appearing as either a suitable restriction of a surjective
Jordan $\ast$-isomorphism, or of the adjoint of such:
\begin{enumerate}
\item Affine bijections from the normal states of $\M_2$ onto the normal
states of $\M_1$.
\item $\textrm{Weak}^{\ast}$ continuous affine bijections from the state space
of $\M_2$ onto that of $\M_1$.
\item Linear identity-preserving isometries from $\M_1$ onto $\M_2$.
\item Linear identity-preserving order-isomorphisms from $\M_1$ onto $\M_2$.
\end{enumerate}
If indeed the von Neumann algebra $\M_1$ contains no direct summand of type
$\textbf{I}_2$, then the set of all orthoisomorphisms from $\p_{\M_1}$ onto
$\p_{\M_2}$ is also in one to one correspondence with the set of surjective
Jordan $\ast$-isomorphisms.
\end{theorem}

\begin{proof}
The one direction of statement (1) follows from the lemma. The other from
a result of Kadison cf. \cite[3.2.8]{BRo}. Statement (2) may be found on p
264 of \cite{Sto} (see the discussion following Corollary 5.8). Statements
(3) and (4) are contained in \cite[3.2.3]{BRo}. The one direction of the
final statement may be deduced from \cite[3.7 \& 4.2]{L2} and the other from
\cite{D} (see p 83).
\end{proof}

Yet another possible description can be had from investigating
those pure state bijections which correspond to  Jordan
$\ast$-isomorphisms between the relevant algebras. From a physical
point of view the importance of such a description is due to the
fact that it demonstrates that there is enough information
internally encoded in the pure states of a $C^\ast$-algebra
(respectively von Neumann algebra) to either guarantee its
physical equivalence to another algebra, or alternatively to be
able to identify a pure state transformation as some sort of
Wigner symmetry. The most impressive result in this regard seems
to be a result of Shultz \cite{Sh} which serves to characterise
those pure state bijections which correspond to surjective
$\ast$-isomorphisms. We review the relevant results before
indicating possible extensions and improvements.

Following Shultz \cite{Sh} we define the transition probability
for pure states as follows:

\begin{definition}
Let $\mathcal{A}$ be a $C^\ast$-algebra and let $\omega_0$ and
$\omega_1$ be pure states. We call $\omega_0$ and $\omega_1$
unitarily equivalent if we can find a unitary element $U \in
\mathcal{A}$ such that $\omega_1(U^\ast(\cdot)U)=
\omega_0(\cdot)$, and unitarily inequivalent if this is not
possible. If $\omega_0$ and $\omega_1$ are unitarily inequivalent
we define the transition probability from $\omega_0$ to $\omega_1$
to be zero. If $\omega_1$ and $\omega_1$ are unitarily equivalent,
there exists an irreducible representation $(\pi, \mathfrak{h})$
and unit vectors $x_0, x_1 \in \mathfrak{h}$ such that $\omega_0 =
\omega_{x_0} \circ \pi$ and $\omega_1 = \omega_{x_1} \circ \pi$
\cite[10.2.3 \& 10.2.6]{KRi}. In this case the transition
probability from $\omega_0$ to $\omega_1$ is defined to be
$\abs{(x_0, x_1)}^2$. By for example \cite[10.3.7]{KRi} this
covers all possibilities.
\end{definition}

In his analysis of pure state transformations Shultz used the
notions of transition probabilities and orientation to
characterise those pure state bijections which correspond to
$\ast$-isomorphisms. The concept of orientation of the set of pure
states is based on the fact that if two distinct pure states are
unitarily equivalent the face of the state space they generate is
a 3-ball, with the face being a line segment if they are
inequivalent. Given this fact the notion of orientation is nothing
but an extension of the fact each 3-ball can effectively be
oriented in one of two ways using either a right-hand or left-hand
set of axes. For further details the reader is referred to the
paper of Shultz \cite{Sh}.

With regard to the following result, recall that the bidual
$\mathcal{A}^{\ast\ast}$ of a $C^\ast$-algebra $\mathcal{A}$ may
be identified with the double commutant of its universal
representation \cite[10.1.1]{KRi}. With this in mind the atomic
part of $\mathcal{A}^{\ast\ast}$ is then the direct sum of all the
direct summands of $\mathcal{A}^{\ast\ast}$ which are type I
factors. However unitarily equivalent pure states of $\mathcal{A}$
generate equivalent irreducible representations \cite[10.2.3 \&
10.2.6]{KRi} and hence with regard to the GNS construction
correspond to $\ast$-isomorphic copies of type I factors in the
direct sum representation of $\mathcal{A}^{\ast\ast}$. If among
these type I factors one ``removes all duplicates'', the result is
what is called the \emph{reduced atomic representation}. This will
be discussed in further detail later on.

\begin{theorem}[\cite{Sh}]
Let $\mathcal{A}$ and $\mathcal{B}$ be (not necessarily unital)
$C^\ast$-algebras. A bijection $\varphi^\sharp$ from the set
$\mathcal{P}_{\mathcal B}$ of pure states of $\mathcal{B}$ onto
$\mathcal{P}_{\mathcal A}$ is induced by a $\ast$-isomorphism
$\varphi$ from the atomic part of $\mathcal{A}^{\ast\ast}$ onto
the atomic part of $\mathcal{B}^{\ast\ast}$ in the sense that
$\varphi^\sharp$ appears as a restriction of the adjoint of
$\varphi$ if and only if $\varphi^\sharp$ preserves both
transition probabilities and orientation. Moreover a bijection
$\varphi^\sharp$ from $\mathcal{P}_{\mathcal B} \cup \{0\}$ onto
$\mathcal{P}_{\mathcal A} \cup \{0\}$ with $\varphi^\sharp(0) =
(0)$ is induced by a $\ast$-isomorphism $\varphi$ from
$\mathcal{A}$ onto $\mathcal{B}$ if and only if $\varphi^\sharp$
is a uniform $\sigma(\mathcal{B}^\ast,
\mathcal{B})$-$\sigma(\mathcal{A}^\ast, \mathcal{A})$
homeomorphism which preserves transition probabilities and
orientation.
\end{theorem}

We will see that if in the above result one deletes the condition
regarding the preservation of orientation the result is a
corresponding relationship between pure state bijections and
Jordan $\ast$-isomorphisms between the relevant algebras. More
importantly we will indicate how one may use the result of Dye
\cite{D} listed in Theorem 1.3 to show that if no irreducible
representation of $\mathcal{A}$ is of the form $M_2(\mathbb{C})$,
then \emph{significantly less} information is needed to get the
same correspondences. In particular in this case we don't need to
know that all transition probabilities are preserved, but rather
only the orthogonal ones. This complements and extends similar
results by Cassinelli, et al (\cite{CdVLL}) who considered only
the case $\mathcal{A} = B(\mathfrak{h})$.

More generally we will precisely describe those transformations
between the pure state spaces of two $C^*$-algebras which in a
canonical way correspond to linear mappings on the underlying
algebras (with no assumptions of bijectivity). To achieve such a
description we also investigate the relationship between
transformations which behave well with regard to unitary
equivalence and orthogonality of pure states. These results of
course extend the work of Shultz \cite{Sh} who considered only the
case of bijections, and focused on the preservation of transition
probabilities rather than orthogonality. We point out that the
results presented in this paper complement the characterisations
of linear maps on $C^*$-algebras with pure state preserving
adjoints given in \cite{Sto1} and \cite{LMa} (\cite[Theorem
5.6]{Sto1} gives a local version of this result and \cite[Theorem
5]{LMa} a global version). In the case of $C(K)$ spaces ($K$
compact Hausdorff) a linear operator $\varphi: C(K) \rightarrow
C(S)$ is called a composition operator if it is induced by some
transformation $T:S \rightarrow K$ in the sense that $\varphi(f) =
f \circ T$ for each $f \in C(K)$. Keeping in mind that up to a
homeomorphism $K$ and $S$ are effectively just the pure state
spaces of $C(K)$ and $C(S)$, these cycles of results therefore
effectively lay the groundwork for a noncommutative theory of
composition operators on $C^*$-algebras.

\section{Orthogonality of pure states}

We analyse the concept of orthogonality of pure states, before
proceeding to investigate transformations between pure state
spaces that behave well with regard to orthogonality. Ultimately
we will see that for such transformations it is precisely their
behaviour with regard to orthogonality together with some
continuity restrictions that determine whether they are induced by
adjoints of linear maps on the underlying algebras or not.

\begin{definition}
Let $\mathcal{A}$ be a $C^\ast$-algebra. Two states $\omega_0$ and
$\omega_1$ are called \emph{orthogonal} if $\|\omega_0 -
\omega_1\| = \|\omega_0\| + \|\omega_1\|$, ie. $\|\omega_0 -
\omega_1\| = 2$.
\end{definition}

\begin{remark}
We note that the concept of orthogonality used here is that
presented in the book of Conway (\cite{Co}), and is slightly
different to the concept of orthogonality presented in the book of
Bratteli and Robinson (\cite{BRo}) being in some sense more
geometrical in flavour. To illustrate the point note that for any
two states $\omega_0$ and $\omega_1$ the fact that all positive
functionals assume their norm at $\I \in \mathcal{A}$ ensures that
$\|\omega_0 + \omega_1\| = \omega_0(\I) + \omega_1(\I) =
\|\omega_0\| + \|\omega_1\|$. Hence for states the above
definition boils down to a concept of orthogonality of $\omega_0$
and $\omega_1$ based on the fact that the vectors $\omega_0 +
\omega_1$ and $\omega_0 - \omega_1$ have the same length.
\end{remark}

The so-called \emph{reduced atomic representation} of a
$C^\ast$-algebra will prove to be an important tool in
establishing the main theorem of this and the next section. For
the sake of clarity of exposition we therefore proceed to review
the most important structural information regarding this
representation. For further details the reader is referred to for
example (\cite{KRi}).

\begin{remark} \label{rar}
If $(\pi, \mathfrak{h})$ is the reduced atomic representation of
the $C^\ast$-algebra $\mathcal{A}$, then $\pi$ is faithful and for
some maximal set $\{(\pi_a, \mathfrak{h}_a)\}_{\mathbb{A}}$ of
pairwise inequivalent irreducible representations, the
representation $\pi(\mathcal{A}) \subset B(\mathfrak{h})$ is of
the form $\pi = \oplus_{a \in \mathbb{A}} \pi_a$ with
$\mathfrak{h} = \oplus_{a \in \mathbb{A}} \mathfrak{h}_a$ and
$\pi(\mathcal{A})'' = \oplus_{a \in \mathbb{A}}
B(\mathfrak{h}_a)$. (See \cite[10.3.10]{KRi} and the discussion
preceding it.) From for example \cite[10.2.3 \& 10.2.5]{KRi} it
is clear that the pure states of $\mathcal{A}$ correspond exactly
to the norm-one vectors of $\oplus_{a \in \mathbb{A}}
\mathfrak{h}_a$ of the form $x = \oplus_{a \in \mathbb{A}} x_a$
where $x_a = 0$ for all but one $a \in \mathbb{A}$. Now given a
vector $x_b$ in $\mathfrak{h}_b$ we will write $\hat{x}_b$ for the
vector in $\mathfrak{h} = \oplus_{a \in \mathbb{A}}
\mathfrak{h}_a$ with the $b^{\mathrm{th}}$ coordinate precisely
$x_b$ and all other coordinates zero. Given pure states $\omega_0$
and $\omega_1$ corresponding to say $\hat{x}_{a_0}$ and
$\hat{z}_{a_1}$, it is clear from \cite[10.3.7]{KRi} that
$\omega_0$ and $\omega_1$ are disjoint if and only if $a_0 \neq
a_1$ (ie. if and only if $\omega_0$ and $\omega_1$ are
inequivalent). Finally considering Definition 2.4 in the context
of the above construction, it is clear that the transition
probability of two pure states $\omega_0$ and $\omega_1$ of
$\mathcal{A}$ corresponding to say $\hat{x}_{a_0}$ and
$\hat{z}_{a_1}$ is precisely $\abs{(\hat{x}_{a_0},
\hat{z}_{a_1})}^2$.
\end{remark}

We need the following facts regarding orthogonal pure states.
Although probably known, the author is not aware of a concise and
explicit statement of these facts, and therefore for the sake of
completeness has elected to reproduce the proof locally.

\begin{proposition} \label{orth}
Let $\mathcal{A}$ be a $C^\ast$-algebra and $\omega_0$ and
$\omega_1$ pure states of $\mathcal{A}$. Then the following are
equivalent:
\begin{enumerate}
\item $\omega_0$ and $\omega_1$ are orthogonal;
\item $\omega_0$ and $\omega_1$ are either disjoint or there
exists an irreducible representation $(\pi_0, \mathfrak{h}_0)$ of
$\mathcal{A}$ and orthonormal vectors $x_0, x_1 \in
\mathfrak{h}_0$ such that $$\omega_0(\cdot) =
\omega_{x_0}\circ\pi_0(\cdot) \quad \omega_1(\cdot) =
\omega_{x_1}\circ\pi_0(\cdot);$$
\item if $(\pi, \mathfrak{h})$ is the reduced atomic representation
of $\mathcal{A}$, we may find orthonormal vectors $x_0, x_1 \in
\mathfrak{h}$ such that $$\omega_0(\cdot) =
\omega_{x_0}\circ\pi(\cdot) \quad \omega_1(\cdot) =
\omega_{x_1}\circ\pi(\cdot);$$
\item if $(\pi, \mathfrak{h})$ is the reduced atomic representation
of $\mathcal{A}$ and if by ${\widetilde{\omega}}_0$ and
${\widetilde{\omega}}_1$ we denote the unique normal extensions of
$\omega_0 \circ \pi^{-1}$ and $\omega_1 \circ \pi^{-1}$ to all of
${\pi(\mathcal{A})}''$ \cite[Lemma 11]{LMa}, there exists an
orthogonal projection $E \in {\pi(\mathcal{A})}''$ such that
$${\widetilde{\omega}}_0(E\pi(\cdot)E) = \omega_0, \quad
{\widetilde{\omega}}_0((\I - E)\pi(\cdot)(\I - E)) = 0$$ and
$${\widetilde{\omega}}_1(E\pi(\cdot)E) = 0, \quad
{\widetilde{\omega}}_1((\I - E)\pi(\cdot)(\I - E)) = \omega_1$$
\end{enumerate}
\end{proposition}

\begin{proof}
Without loss of generality we may identify $\mathcal{A}$ with its
reduced atomic representation.

\underline{$(2) \Leftrightarrow (3)$}: $\quad$ In the light of the
information garnered from Remark \ref{rar} this is a
straightforward exercise.

\underline{$(3) \Leftrightarrow (4)$}: $\quad$ Let $\omega_0$ and
$\omega_1$ be pure states respectively corresponding to norm-one
vectors $\hat{x}_{a_0}$ and $\hat{z}_{a_1}$ in the sense described
in Remark \ref{rar}.

First suppose that (4) holds, that is there exists a projection $E
\in {\mathcal{A}}''$ with $$0 =\omega_0(\I - E) =
\omega_{\hat{x}_{a_0}}(\I - E) = \|(\I - E)\hat{x}_{a_0}\|^2$$ and
$$0 =\omega_1(E) = \omega_{\hat{z}_{a_1}}(E) =
\|E\hat{z}_{a_1}\|^2.$$ This can clearly only be the case if
$\hat{x}_{a_0} \perp \hat{z}_{a_1}$.

Conversely if $\hat{x}_{a_0} \perp \hat{z}_{a_1}$, we may set $E =
E_0$ where $E_0$ is the orthogonal projection of $\mathfrak{h} =
\oplus_{a \in \mathbb{A}} \mathfrak{h}_a$ onto
$\mathrm{span}\{\hat{x}_{a_0}\}$. It is not difficult to see that
$E$ then belongs to $\mathcal{A}'' = \oplus_{a \in \mathbb{A}}
B(\mathfrak{h}_a)$. Moreover for each $A \in \mathcal{A}$ we have
that $$\omega_0(EAE) = (AE\hat{x}_{a_0}, E\hat{x}_{a_0}) =
(A\hat{x}_{a_0}, \hat{x}_{a_0}) = \omega_0(A)$$ and
$$\omega_1(EAE) = (AE\hat{z}_{a_1}, E\hat{z}_{a_1}) = 0.$$
Similarly $$\omega_0((\I - E)\cdot(\I - E)) = 0 \quad \mbox{and}
\quad \omega_1((\I - E)\cdot(\I - E)) = \omega_1(\cdot).$$

\underline{$(4) \Rightarrow (1)$}: $\quad$ Let $\omega_0,
\omega_1$ be pure states which satisfy condition (4) for some
orthogonal projection $E \in \mathcal{A}''$. Since by assumption
$\omega_0$ and $\omega_1$ are normal, each of $\omega_0, \omega_1$
and $\omega_0 - \omega_1$ extends uniquely and without change of
norm to $\mathcal{A}''$ \cite[10.1.11]{KRi}. Now since $2E - \I$
is self-adjoint with $(2E - \I)^2 = \I$, we clearly have $\|2E -
\I\| = 1$. But then \begin{eqnarray*} \|\omega_0 - \omega_1\|
&\geq& (\omega_0 - \omega_1)(2E - \I) = (\omega_0 - \omega_1)(E -
(\I - E))\\ &=& \omega_0(E) + \omega_1(\I - E) = \omega_0(\I) +
\omega_1(\I) = \|\omega_0\| + \|\omega_1\|
\end{eqnarray*}
by condition (4). Since in general $\|\omega_0 - \omega_1\| \leq
\|\omega_0\| + \|\omega_1\|$, we conclude that $\|\omega_0 -
\omega_1\| = \|\omega_0\| + \|\omega_1\|$.

\underline{$(1) \Rightarrow (2)$}: $\quad$ We verify this
implication for vector states of $B(\mathfrak{k})$ ($\mathfrak{k}$
an arbitrary Hilbert space). By Remark \ref{rar} it is clear that
this will suffice to establish the implication for the case of
equivalent pure states of $\mathcal{A}$. Thus let $x, y \in
\mathfrak{k}$ be unit vectors, let $\omega_0 = \omega_x$ and
$\omega_1 = \omega_y$, and assume that $\|\omega_0 - \omega_1\| =
2$. Now by replacing $y$ with $\widetilde{y} = e^{i\theta}y$ if
necessary where $(x, y) = \abs{(x, y)}e^{i\theta}$, we may assume
without loss of generality that $(x, y) = \abs{(x, y)}$. Since
vector states are clearly $\sigma$-weakly continuous, $\omega_0$
and $\omega_1$ belong to the predual of $B(\mathfrak{k})$, and
hence the Hahn-Banach theorem assures us that there exists $A \in
(B(\mathfrak{k})_\ast)^\ast = B(\mathfrak{k})$ with $$\|A\| = 1
\quad \mbox{and} \quad \omega_0(A) - \omega_1(A) = \|\omega_0 -
\omega_1\| = 2.$$ Now since $\abs{\omega_0(A)} \leq \|\omega_0\|
\|A\| = 1$ and $\abs{\omega_1(A)} \leq 1$, this can only be the
case if $$\omega_0(A) = 1 \quad \mbox{and} \quad \omega_1(A) =
-1.$$ On bringing the Cauchy-Schwarz inequality into play we may
now conclude that $$1 = \omega_0(A) = (Ax, x) \leq \|Ax\| \|x\|
\leq \|A\| \|x\|^2 = 1$$ and hence that $(Ax, x) = \|Ax\| = \|x\|
= 1$. Consequently $$\|Ax - x\|^2 = \|Ax\|^2 - 2\Re(Ax, x) +
\|x\|^2 = 0,$$ ie. $Ax = x$. Similarly $Ay = -y$. But then $$\|x -
y\| = \|A(x + y)\| \leq \|x + y\| = \|A(x - y)\| \leq \|x - y\|,$$
that is $\|x - y\| = \|x + y\|$. Squaring both sides and
cancelling like terms reveals that $-2\Re(x, y) = 2\Re(x, y)$.
Since by assumption $(x, y) = \abs{(x, y)}$, we conclude that $(x,
y) = \Re(x, y) = 0$ as required.
\end{proof}

\begin{definition}
Let $\mathcal{A, B}$ be $C^*$-algebras. A transformation
$\varphi^\sharp$ from a subset of the pure states of $\mathcal{B}$
into the set of pure states of $\mathcal{A}$, is called
\emph{orthogonal} if for any two pure states $\omega_0, \omega_1
\in \mathrm{dom}(\varphi^{\sharp})$, we have that
$\varphi^{\sharp}(\omega_0)$ and $\varphi^{\sharp}(\omega_1)$ are
orthogonal whenever $\omega_0$ and $\omega_1$ are orthogonal.

If on the other hand $\omega_0$ and $\omega_1$ are orthogonal
whenever $\varphi^{\sharp}(\omega_0)$ and
$\varphi^{\sharp}(\omega_1)$ are orthogonal, we call
$\varphi^\sharp$ co-orthogonal.

If $\varphi^\sharp$ is both orthogonal and co-orthogonal, it will
be deemed to be bi-orthogonal.
\end{definition}

\begin{definition}
Let $\mathcal{A, B}$ be $C^\ast$-algebras and let $\varphi^\sharp$
be a transformation from the set of pure states of $\mathcal{B}$
into the set of pure states of $\mathcal{A}$. For a pure state
$\omega$ the set of all other pure states unitarily equivalent to
$\omega$ will be denoted by $[\omega]$. We call $\varphi^\sharp$
\emph{fibre-preserving} if $\varphi^\sharp([\omega]) \subset
[\varphi^\sharp(\omega)]$ for each pure state $\omega$ of
$\mathcal{B}$.

Given a map $\varphi^\sharp$ from $\mathcal{P}_{\mathcal{B}}$ into
$\mathcal{P}_{\mathcal{A}}$, we say that a property holds
\emph{locally} for $\varphi^\sharp$ if it holds for the
restriction of $\varphi^\sharp$ to each equivalence class
$[\omega]$ of $\mathcal{P}_{\mathcal{B}}$. For example we call
$\varphi^\sharp$ locally injective if for each $\omega \in
\mathcal{P}_{\mathcal{B}}$ $\varphi^\sharp$ acts injectively on
$[\omega]$, locally orthogonal if for each $\omega \in
\mathcal{P}_{\mathcal{B}},$ $\varphi^\sharp|_{[\omega]}$ is
orthogonal, etc.
\end{definition}

\begin{definition}
Let $\mathcal{A, B}$ be $C^\ast$-algebras and let $\varphi^\sharp$
be a fibre-preserving transformation from
$\mathcal{P}_{\mathcal{B}}$ into $\mathcal{P}_{\mathcal{A}}$. For
a subset $\mathcal{V}$ of say $\mathcal{P}_{\mathcal{A}}$, we
write $\mathcal{V}^\perp$ for the set $$\mathcal{V}^\perp =\{
\omega \in \mathcal{P}_{\mathcal{A}} | \omega \perp \rho, \rho \in
\mathcal{V}\}.$$ We call the range of $\varphi^{\sharp}$
\emph{locally solid} if $\varphi^\sharp([\omega]) =
\varphi^\sharp([\omega])^{\perp \perp}$ for each $\omega \in
\mathcal{P}_{\mathcal{B}}$.
\end{definition}

\begin{remark}
Let $\mathcal{A}$ be a commutative $C^\ast$-algebra. In this
setting it follows from for example \cite[4.4.1]{KRi} that pure
states $\omega_0$ and $\omega_1$ of $\mathcal{A}$ are unitarily
equivalent if and only if $\omega_0 = \omega_1$. Thus
``equivalence classes'' of pure states are just singletons. Since
in general two pure states are inequivalent if and only if they
are disjoint \cite[10.2.3 \& 10.3.7]{KRi}, we may apply
Proposition \ref{orth}(2) to the above fact to conclude that in
the commutative setting pure states $\omega_0$ and $\omega_1$ are
orthogonal if and only if $\omega_0 \neq \omega_1$. Consequently
any transformation between the pure state spaces of commutative
$C^\ast$-algebras is automatically both co-orthogonal and
fibre-preserving, with a transformation being orthogonal precisely
when it is injective. The orthogonality of distinct pure states
(point evaluations) in the case $\mathcal{A} = C(K)$, may also be
deduced from Urysohn's lemma. So in a very real sense the
requirement of orthogonality in the non-commutative setting
compensates for the lack of a suitable ``non-commutative''
Urysohn's lemma.
\end{remark}

Our primary task is of course to identify those transformations
$\varphi^\sharp : \mathcal{P}_{\mathcal{B}} \rightarrow
\mathcal{P}_{\mathcal{A}}$ that are induced by linear maps on the
underlying algebras. Taking our cue from the commutative case
elucidated above, a good place to start seems to be among the
co-orthogonal fibre-preserving tranformations. We therefore
proceed to investigate the relationship between these two
properties.

\begin{lemma}\label{coorth}
Let $\mathcal{A, B}$ be $C^*$-algebras and $\varphi^\sharp$
a transformation from $\mathcal{P}_{\mathcal{B}}$ into
$\mathcal{P}_{\mathcal{A}}$. Consider the following statements:
\begin{enumerate}
\item $\varphi^\sharp$ is co-orthogonal.
\item $\varphi^\sharp$ is locally co-orthogonal.
\item $\varphi^\sharp$ is fibre-preserving.
\end{enumerate}
The implications $(1) \Leftrightarrow (2) \Rightarrow (3)$ hold in
general.
\end{lemma}

\begin{proof} The implication $(1) \Rightarrow (2)$ is of course
trivial. We therefore proceed to verify that $(2) \Rightarrow
(3)$. Firstly assume that both $\mathcal{A, B}$ are reduced
atomically represented and that $\varphi^\sharp$ is locally
co-orthogonal. Now let $\omega_0, \omega_1 \in
\mathcal{P}_{\mathcal B}$ be unitarily equivalent pure states. We
show that $\varphi^\sharp(\omega_0)$ and
$\varphi^\sharp(\omega_1)$ are then necessarily also unitarily
equivalent. In the notation of Remark \ref{rar} $\omega_0$ and
$\omega_1$ are of the form $\omega_0 = \omega_{\hat{x}_a}$ and
$\omega_1 = \omega_{\hat{y}_a}$ where for some fixed $a \in
\mathbb{B}$ both $x_a$ and $y_a$ are unit vectors in
$\mathfrak{k}_a$. If now $\omega_0$ and $\omega_1$ are not
orthogonal, then by the local co-orthogonality of
$\varphi^\sharp$, neither are $\varphi^\sharp(\omega_0)$ and
$\varphi^\sharp(\omega_1)$. Hence in this case
$\varphi^\sharp(\omega_0)$ and $\varphi^\sharp(\omega_1)$ are
equivalent. If however $\omega_0$ and $\omega_1$ are indeed
orthogonal but unitarily equivalent pure states (ie. $x_a \perp
y_a$), then $\omega_2 = \omega_{\hat{z}_a}$ where $z_a =
\frac{1}{\sqrt{2}}(x_a + y_a) \in \mathfrak{k}_a$ is a pure state
of $\mathcal{B}$ which by Proposition \ref{orth} is orthogonal to
neither $\omega_0$, nor $\omega_1$. Again the local
co-orthogonality of $\varphi^\sharp$ now ensures that
$\varphi^\sharp(\omega_2)$ is orthogonal to neither
$\varphi^\sharp(\omega_0)$ nor $\varphi^\sharp(\omega_1)$. Thus
both $\varphi^\sharp(\omega_0)$ and $\varphi^\sharp(\omega_1)$ are
unitarily equivalent to $\varphi^\sharp(\omega_2)$, and hence
equivalent to each other.

It remains to show that $(2) \Rightarrow (1)$. To this end let
$\varphi^\sharp$ be locally co-orthogonal and let $\omega_0$ and
$\omega_1$ be pure states of $\mathcal{B}$ with
$\varphi^\sharp(\omega_0)$ and $\varphi^\sharp(\omega_1)$
orthogonal. If $\omega_0$ is unitarily equivalent to $\omega_1$,
then by the local co-orthogonality of $\varphi^\sharp$, $\omega_0$
and $\omega_1$ must also be orthogonal. If $\omega_0$ is not
unitarily equivalent to $\omega_1$, then we necessarily already
have that $\omega_0$ and $\omega_1$ are orthogonal. Thus
$\varphi^\sharp$ is co-orthogonal.
\end{proof}

Again as in the commutative case there is a clear link between
orthogonality and injectivity.

\begin{lemma} \label{inj} Let $\mathcal{A, B}$ and $\varphi^\sharp$ be as in
the previous lemma. If $\varphi^\sharp$ is bi-orthogonal
(alternatively locally bi-orthogonal), then it is injective
(locally injective).
\end{lemma}

\begin{proof}
Due to the similarity of the proofs we prove only the first claim.
Hence assume that both $\mathcal{A, B}$ are reduced atomically
represented and that $\varphi^\sharp$ is bi-orthogonal. Let
$\omega_0, \omega_1 \in \mathcal{P}_{\mathcal B}$ be two distinct
pure states. If $\omega_0$ and $\omega_1$ are orthogonal, then by
the bi-orthogonality of $\varphi^\sharp$, so are
$\varphi^\sharp(\omega_0)$ and $\varphi^\sharp(\omega_1)$. Thus in
this case $\varphi^\sharp(\omega_0)$ and
$\varphi^\sharp(\omega_1)$ are clearly distinct. If now $\omega_0$
and $\omega_1$ are not orthogonal, they must of course necessarily
be unitarily equivalent. In the notation of Remark \ref{rar}
$\omega_0$ and $\omega_1$ are of the form $\omega_0 =
\omega_{\hat{x}_a}$ and $\omega_1 = \omega_{\hat{y}_a}$ where for
some fixed $a \in \mathbb{B}$ both $x_a$ and $y_a$ are unit
vectors in $\mathfrak{k}_a$. Now by Proposition \ref{orth} the two
vectors $x_a$ and $y_a$ are distinct but not orthogonal. Thus in
the two-dimensional subspace spanned by these vectors we can find
a third unit vector $z_a$ which is orthogonal to $x_a$ but not to
$y_a$. By Proposition \ref{orth} this unit vector corresponds to
pure state $\omega_2 \in \mathcal{P}_{\mathcal{B}}$ which is
unitarily equivalent to both $\omega_0$ and $\omega_1$, orthogonal
to $\omega_0$, but not to $\omega_1$. Since $\varphi^\sharp$ is
(locally) bi-orthogonal, it follows that
$\varphi^\sharp(\omega_2)$ is orthogonal to
$\varphi^\sharp(\omega_0)$, but not to $\varphi^\sharp(\omega_1)$.
This can clearly only be the case if $\varphi^\sharp(\omega_0)$
and $\varphi^\sharp(\omega_1)$ are distinct. Thus $\varphi^\sharp$
is injective.
\end{proof}

We have already proposed co-orthogonality and the concomitant
preservation of fibres as properties that may help to identify
those pure state transformations that come from linear mappings on
the underlying algebras. In addition to these some local
properties will no doubt also be needed. However since in the
commutative case (unitary) equivalence classes of pure states are
just singletons, it is not so easy to use this case to formulate a
conjecture regarding the necessary local properties. We therefore
proceed to investigate the case of $B(\mathfrak{h})$ in order to
get some idea of what may be needed. Classically Wigner used
preservation of transition probabilities to identify those vector
state bijections which come from linear maps on $B(\mathfrak{h})$
(\cite{W}; cf.\cite[Theorem 1]{Sh}). The explicit use of
orthogonality to achieve the same end, seems to be a result due to
Cassinelli, de Vito, et al \cite{CdVLL}. We sketch a proof of
their result before extending both these results to the most
general (non-bijective) case possible. To prove the bijective case
we need the following lemma. Although well-known, the author has
not been able to find an explicit statement of the anti-isomorphic
case, and hence has once again elected to reproduce the proof
locally.

\begin{lemma}
Let $\varphi$ be a $\ast$-isomorphism or $\ast$-antiisomorphism
from $B(\mathfrak{h})$ onto $B(\mathfrak{k})$. Then
$$Tr_{\mathfrak{h}} = Tr_{\mathfrak{k}} \circ \varphi.$$
\end{lemma}

\begin{proof}
Without loss of generality assume that $\mathrm{dim}(\mathfrak{h})
= \infty$. First note that with $\varphi$ as above we may find a
unitary operator $U: \mathfrak{k} \rightarrow \mathfrak{h}$ such
that either
\begin{equation} \label{eq:c}
\varphi(A) = U^\ast A U \quad \mbox{for all} \quad A \in
B(\mathfrak{h})
\end{equation}
or
\begin{equation} \label{eq:d}
\varphi(A) = U^\ast c^\ast A^\ast c U \quad \mbox{for all} \quad A
\in B(\mathfrak{h})
\end{equation}
where $c$ is the anti-unitary operator on $\mathfrak{h}$ induced
by complex conjugation. To see this one may for example suitably
adapt the proof of \cite[Example 3.2.14]{BRo}. (In this regard
note that the map $\sigma_0$ in \cite[Example 3.2.14]{BRo}
corresponds to transposition with respect to some orthonormal base
and hence to all intents of purposes is of the form $\sigma_0(A) =
c^\ast A^\ast c$.

Alternatively one may apply \cite[Theorem 16.B(1)]{LMa} to
conclude that there exist injective partial isometries $U:
\mathfrak{k} \rightarrow \mathfrak{h}$ and $V: \mathfrak{k}
\rightarrow \mathfrak{h}$ such that either
\begin{equation} \label{eq:a}
\varphi(A) = V^\ast A U \quad \mbox{for all} \quad A \in
B(\mathfrak{h})
\end{equation}
or
\begin{equation} \label{eq:b}
\varphi(A) = V^\ast c^\ast A^\ast c U \quad \mbox{for all} \quad A
\in B(\mathfrak{h})
\end{equation}
where $c$ is as before. (The surjectivity of $\varphi$ excludes
part $B(2)$ of \cite[Theorem 16.B]{LMa} as a possibility.) To
see that \cite[Theorem 16.B(1)]{LMa} is indeed applicable we may
combine Lemma 1.2, Theorem 1.5 and \cite[Corollary 20]{LMa} to
conclude that $\varphi^{\ast}$ maps the $\sigma$-weakly continuous
extreme points of the dual ball of $B(\mathfrak{k})$ onto extreme
points of the dual ball of $B(\mathfrak{h})$. Having obtained the
above formulae one may then note that $\varphi(\I) = \varphi(E)$
where $E$ is the range projection of $U$, and then conclude from
the injectivity of $\varphi$ that $\I = E$. Thus $U$ must in fact
be surjective, and hence a unitary. It is then a simple matter to
see that since $\I = \varphi(\I)$, we must then have that $V^\ast
= V^{\ast}UU^\ast = \varphi(\I)U^\ast = {\I}U^\ast = U^\ast$.

Now for any set of orthonormal bases $\{x_\lambda\}_\Lambda$ and
$\{y_\lambda\}_\Lambda$ of $\mathfrak{k}$,
$\{{\widetilde{x}}_\lambda\}_\Lambda = \{cU(x_\lambda)\}_\Lambda$
and $\{{\widetilde{y}}_\lambda\}_\Lambda =
\{cU(y_\lambda)\}_\Lambda$ are orthonormal bases of
$\mathfrak{h}$. If therefore $\varphi$ is of the form
(\ref{eq:d}), then
\begin{eqnarray*}
\sum_{\lambda \in \Lambda} (A{\widetilde{x}}_\lambda,
{\widetilde{y}}_\lambda) &=& \sum_{\lambda \in \Lambda}
(AcUx_\lambda, cUy_\lambda)\\ &=& \sum_{\lambda \in \Lambda}
(Uy_\lambda, c^\ast(AcUx_\lambda))\\ &=& \sum_{\lambda \in
\Lambda} (y_\lambda, U^{\ast}c^{\ast}AcUx_\lambda)\\ &=&
\sum_{\lambda \in \Lambda} (y_\lambda, \varphi(A^\ast)x_\lambda)\\
&=& \sum_{\lambda \in \Lambda} (\varphi(A)y_\lambda, x_\lambda)
\end{eqnarray*}
for each $A \in B(\mathfrak{h})$. In particular it then follows
that $\varphi(A)$ is trace-class whenever $A$ is trace class, and
that $$Tr_{\mathfrak{h}}(A) = \sum_{\lambda \in \Lambda}
(A{\widetilde{x}}_\lambda, {\widetilde{x}}_\lambda) =
\sum_{\lambda \in \Lambda} (\varphi(A)x_\lambda, x_\lambda) =
Tr_{\mathfrak{k}}(\varphi(A))$$ for each trace-class element of
$B(\mathfrak{h})$. A similar conclusion holds if $\varphi$ is of
the form (\ref{eq:c}).
\end{proof}

The one-dimensional subspaces of $\mathfrak{h}$ are clearly in a
one-one correspondence with the vector states of
$B(\mathfrak{h})$. More precisely if $x \in \mathfrak{h}$ is a
unit vector and $E_x$ the minimal projection onto the ray spanned
by $x$, then $\omega_x = Tr_{\mathfrak{h}}(E_x\cdot E_x)$. Thus
with reference to Definition 1.4 and the discussion preceding
Theorem 1.3, the following result is the first step towards
showing that on condition we exclude the case of
$M_2(\mathbb{C})$, we don't need to know that all transition
probabilities are preserved in order to identify a pure state
transformation as a Wigner symmetry. As noted earlier this result
is of course known (see for example Cassinelli, de Vito, Lahti and
Levrero \cite{CdVLL}). We therefore content ourselves with merely
sketching how this result may be deduced from the previous lemma
by means of Dye's result \cite[p. 83]{D}. We hasten to add that
the dimensional restriction \emph{can not be removed} as the
implication may in fact fail in the two-dimensional case
(\cite{U}; cf \cite[Example 4.1]{CdVLL}).

\begin{theorem}[\cite{CdVLL}] \label{bijbh}
Let $\varphi^\sharp$ be a bijection from the set of vector states
of $B(\mathfrak{k})$ onto the set of vector states of
$B(\mathfrak{h})$. If either $\mathrm{dim}(\mathfrak{h}) \neq 2$
or $\mathrm{dim}(\mathfrak{k}) \neq 2$, then $\varphi^\sharp$ is
\emph{bi-orthogonal} in the sense that for any norm-one vectors
$x, y \in \mathfrak{k}$ and $\widetilde{x}, \widetilde{y} \in
\mathfrak{h}$ with $$\varphi^\sharp(\omega_x) =
\omega_{\widetilde{x}} \quad \mbox{and} \quad
\varphi^\sharp(\omega_y) = \omega_{\widetilde{y}},$$ we always
have that $$x \perp y \quad \mbox{if and only if} \quad
\widetilde{x} \perp \widetilde{y}$$ if and only if there exists
either a $\ast$-isomorphism or a $\ast$-antiisomorphism $\varphi$
from $B(\mathfrak{h})$ onto $B(\mathfrak{k})$ such that
$$\varphi^\sharp(\omega_x) = \omega_x \circ \varphi \quad
\mbox{for each norm-one} \quad x \in \mathfrak{k}.$$
\end{theorem}

\begin{proof}
Let $\varphi^\sharp$ be a bijection between the respective sets of
vector states satisfying the stated hypothesis and suppose that
$B(\mathfrak{h}) \neq M_2(\mathbb{C})$. We show that
$\varphi^\sharp$ canonically induces a bijection $\varphi$ from
the rank-one orthogonal projections of $B(\mathfrak{h})$ onto
those of $B(\mathfrak{k})$ in a way that ``preserves
orthogonality'', before proceeding to show that in fact $\varphi$
extends to an orthoisomorphism from $\p_{B(\mathfrak{h})}$ onto
$\p_{B(\mathfrak{k})}$. Now for any norm-one vector $x \in
\mathfrak{k}$ we will write $\widetilde{x}$ for the corresponding
norm-one vector in $\mathfrak{h}$ such that
$\varphi^\sharp(\omega_x) = \omega_{\widetilde{x}}$. Given a
norm-one vector $x \in \mathfrak{k}$, we set
$\varphi(E_{\widetilde{x}}) = E_x$. (Here $E_x$ and
$E_{\widetilde{x}}$ are the orthogonal projections onto
$\mathrm{span}\{x\}$ and $\mathrm{span}\{\widetilde{x}\}$
respectively.) It is now an exercise to show that the fact that
$E_{\widetilde{x}}E_{\widetilde{y}} = 0$ if and only if
$\varphi(E_{\widetilde{x}})\varphi(E_{\widetilde{y}}) = 0$ is
inherited from the so-called bi-orthogonality of $\varphi^\sharp$.
To see that $\varphi$ extends to an orthoisomorphism from
$\p_{B(\mathfrak{h})}$ onto $\p_{B(\mathfrak{k})}$ we need only
show that by means of $\varphi^\sharp$ we may identify the closed
linear subspaces of $\mathfrak{h}$ with those of $\mathfrak{k}$ in
a one-to-one way that preserves mutual orthogonality of subspaces.
To this end let $\mathfrak{h}_0$ be a closed linear subspace of
$\mathfrak{h}$ and let $\{{\widetilde{x}}_\lambda\}_\Lambda
\subset \mathfrak{h}$ be a set of unit vectors such that
$\mathfrak{h}_0 =
\overline{\mathrm{span}}\{{\widetilde{x}}_\lambda\}_\Lambda$. Now
select $\{x_\lambda\}_\Lambda \subset \mathfrak{k}$ so that
$\varphi^\sharp(\omega_{x_\lambda}) =
\omega_{\widetilde{x}_\lambda}$ for each $\lambda \in \Lambda$.
Now with $\mathfrak{k}_0 =
\overline{\mathrm{span}}\{{x}_\lambda\}_\Lambda$, set
$\varphi(E_{\mathfrak{h}_0}) = E_{\mathfrak{k}_0}$. (Here
$E_{\mathfrak{h}_0}$ and $E_{\mathfrak{k}_0}$ respectively denote
the orthogonal projections onto $\mathfrak{h}_0$ and
$\mathfrak{k}_0$.) We show that $E_{\mathfrak{k}_0}$ is uniquely
defined and that $\varphi(E_{{\mathfrak{h}_0}^\perp}) =
E_{{\mathfrak{k}_0}^\perp}$.

Let $\{\widetilde{z}_\mu\}$ be an orthonormal base for
${\mathfrak{h}_0}^\perp$ with as before $\{z_\mu\}$ selected so
that $\varphi^\sharp(\omega_{z_\mu}) = \omega_{\widetilde{z}_\mu}$
for each $\mu$. By the hypothesis $\{z_\mu\}$ is again an
orthonormal system with $\{x_\lambda | \lambda \in \Lambda\}
\subset (\{z_\mu\})^\perp$ and hence $\mathfrak{k}_0 \subset
(\{z_\mu\})^\perp$. To achieve our stated objective we only need
to show that in fact $\mathfrak{k}_0 = (\{z_\mu\})^\perp$. (The
uniqueness of $\mathfrak{k}_0$ then follows from the fact that
$\varphi^\sharp$ identifies any set of norm-one vectors generating
$(\{\widetilde{z}_{\mu}\})^\perp$ with a set of norm-one vectors
generating $(\{z_\mu\})^\perp$. By interchanging the roles of
$\mathfrak{h}_0$ and ${\mathfrak{h}_0}^\perp$ it then also follows
from this that $\varphi^\sharp$ identifies
${\mathfrak{h}_0}^\perp$ with ${\mathfrak{k}_0}^\perp$.) Since
$\mathfrak{k}_0 \subset (\{z_\mu\})^\perp$, it suffices to show
that ${\mathfrak{k}_0}^\perp \cap (\{z_\mu\})^\perp = \{0\}$ in
order to see that $\mathfrak{k}_0 = (\{z_\mu\})^\perp$. Now if it
were possible to find a norm-one vector $$y \in
{\mathfrak{k}_0}^\perp \cap (\{z_\mu\})^\perp = \{x_\lambda |
\lambda \in \Lambda\}^\perp \cap \{z_\mu\}^\perp$$ then for
$\widetilde{y} \in \mathfrak{h}$ with $\varphi^\sharp(\omega_y) =
\omega_{\widetilde{y}}$ we would by the hypothesis have to have
$$\widetilde{y} \in  \{\widetilde{x}_\lambda | \lambda \in
\Lambda\}^\perp \cap \{\widetilde{z}_\mu\}^\perp =
{\mathfrak{h}_0}^\perp \cap {\mathfrak{h}_0}^{\perp\perp};$$ a
situation which is clearly impossible. Thus as claimed $\varphi$
canonically extends to an orthoisomorphism.

By a result of Dye \cite[p. 83]{D} $\varphi$ then extends further
to either a $\ast$-isomorphism or $\ast$-antiisomorphism from
$B(\mathfrak{h})$ onto $B(\mathfrak{k})$. This extension will also
be denoted by $\varphi$. It remains to show that $\varphi^\sharp$
appears as a restriction of the adjoint of $\varphi$. To this end
let $x \in \mathfrak{k}$ be a norm-one vector, and $\widetilde{x}$
the related norm-one vector in $\mathfrak{h}$. By construction we
have that $\varphi(E_{\widetilde{x}}) = E_x$. It is an exercise to
conclude that $\omega_x(\cdot) = Tr_{\mathfrak{k}}(E_x \cdot E_x)$
and $\omega_{\widetilde{x}}(\cdot) =
Tr_{\mathfrak{h}}(E_{\widetilde{x}} \cdot E_{\widetilde{x}})$. By
making use of the previous lemma we may now conclude that
\begin{eqnarray*}
\varphi^{\sharp}(\omega_x)(A) &=& \omega_{\widetilde{x}}(A)\\
 &=& Tr_{\mathfrak{h}}(E_{\widetilde{x}}AE_{\widetilde{x}})\\
 &=& Tr_{\mathfrak{h}} \circ \varphi^{-1}
 (\varphi(E_{\widetilde{x}})\varphi(A)\varphi(E_{\widetilde{x}}))\\
 &=& Tr_{\mathfrak{k}}(E_{x}\varphi(A)E_x)\\
 &=& \omega_x(\varphi(A)) \quad \mbox{for each} \quad A \in
 B(\mathfrak{h}).
\end{eqnarray*}

For the converse note that as in the previous lemma, we may show
that $\varphi$ is either of the form
\begin{equation*}
\varphi(A) = U^\ast A U \quad \mbox{for all} \quad A \in
B(\mathfrak{h})
\end{equation*}
or
\begin{equation*}
\varphi(A) = U^\ast c^\ast A^\ast c U \quad \mbox{for all} \quad A
\in B(\mathfrak{h}),
\end{equation*}
where $c$ is as before and $U: \mathfrak{k} \rightarrow
\mathfrak{h}$ is a unitary. Using this description it is a simple
matter to verify the claim regarding the bi-orthogonality of
$\varphi^\sharp$. The result follows.
\end{proof}

Let $\varphi^\sharp$ be a transformation from the set of vector
states of $B(\mathfrak{k})$ into the set of vector states of
$B(\mathfrak{h})$. We recall from the discussion preceding Theorem
\ref{bijbh} that the vector states are all states of the form
$\omega = Tr(E_{\omega}\cdot E_{\omega})$ for some minimal
projection. By analogy with the earlier definitions for pure
states we say that the range of $\varphi^\sharp$ is \emph{locally
solid} if for every minimal projection $E_0 \leq \vee\{E_\rho|
\rho = \varphi^\sharp(\omega), \omega \mbox{ a vector state of
}B(\mathfrak{k})\}$ there exists a vector state $\omega_0$ of
$B(\mathfrak{k})$ such that $\varphi^\sharp(\omega_0) =
Tr_{\mathfrak{h}}(E_0\cdot E_0)$.

At the Hilbert space level the set $\vee\{E_\rho| \rho =
\varphi^\sharp(\omega), \omega \mbox{ a vector state of
}B(\mathfrak{k})\}$ is of course nothing but the closure of
$\mathfrak{h}_0 = \mathrm{span}\{\widetilde{x} \in \mathfrak{h} |
\omega_{\widetilde{x}} = \varphi^\sharp(\omega), \omega \mbox{ a
vector state of }B(\mathfrak{k}) \}.$ In this context local
solidness is then nothing more than the claim that any unit vector
$x \in \overline{\mathfrak{h}_0}$ must in fact belong to
$\mathfrak{h}_0$, that is that $\mathfrak{h}_0^{\perp\perp} =
\mathfrak{h}_0$.

\begin{theorem}[Generalised Wigner Theorem]\label{gwt}
Let $\varphi^\sharp$ be a transformation from the set of vector
states of $B(\mathfrak{k})$ into the set of vector states of
$B(\mathfrak{h})$. Consider the following statements:
\begin{enumerate}
\item There exists a linear map $\varphi : B(\mathfrak{h})
\rightarrow B(\mathfrak{k})$ such that $\varphi^\sharp(\omega_x) =
\omega_x \circ \varphi$ for every vector state of
$B(\mathfrak{k})$.

\item $\varphi^\sharp$ preserves transition probabilities and has
a locally solid range.

\item $\varphi^\sharp$ is bi-orthogonal (in the sense defined
in the previous theorem) and has a locally solid range.
\end{enumerate}
The implications $(1) \Leftrightarrow (2) \Rightarrow (3)$ hold in
general with all three statements being equivalent if
$\mathrm{dim}(\mathfrak{k}) \neq 2$.

Moreover any linear map $\varphi : B(\mathfrak{h})
\rightarrow B(\mathfrak{k})$ which induces a vector state
transformation $\varphi^\sharp$ in the manner described above is
either of the form
$$\varphi(A) = U^\ast A U \quad \mbox{for all} \quad A \in
B(\mathfrak{h})$$
or
$$\varphi(A) = U^\ast c^\ast A^\ast c U \quad \mbox{for all}
\quad A \in B(\mathfrak{h})$$
where $c$ is the anti-unitary operator on $\mathfrak{h}$ induced
by complex conjugation and $U$ a linear isometry from $\mathfrak{k}$
into $\mathfrak{h}$
\end{theorem}

\begin{proof}
\underline{$(1) \Rightarrow (2)$}: $\quad$ Suppose we can find
such a linear map $\varphi: B(\mathfrak{h}) \rightarrow
B(\mathfrak{k})$ with $$\varphi^\sharp(\omega_x) = \omega_x \circ
\varphi$$ for every vector state of $B(\mathfrak{k})$. Then
$\varphi$ is clearly positivity preserving (since $A \in
B(\mathfrak{h})^+$ then implies that $(\varphi(A)x, x) \geq 0$ for
each $x \in \mathfrak{k})$ and hence bounded. It now easily
follows from the hypothesis that the dual $\varphi^*$ will map the
normal states of $B(\mathfrak{k})$ (the norm-closed convex hull of
the vector states \cite[7.1.12 \& 7.1.13]{KRi}) into the normal
states of $B(\mathfrak{h})$. On regarding say $B(\mathfrak{k})_*$
as a subspace of $B(\mathfrak{k})$, this means that
$\varphi^*(B(\mathfrak{k})_*) \subset B(\mathfrak{h})_*$, and
hence that $\varphi^*: B(\mathfrak{k})^* \rightarrow
B(\mathfrak{h})^*$ restricts to a map from $B(\mathfrak{k})_*$
into $B(\mathfrak{h})_*$. It is now an exercise to show that the
dual of this induced map is exactly $\varphi$. Therefore $\varphi$
is in fact a dual operator and hence necessarily weak* to weak*
continuous. We may therefore apply \cite[Lemma 5.4]{Sto1} to see
that there exists a linear isometry $U$ from $\mathfrak{k}$ into
$\mathfrak{h}$ such that either
\begin{equation} \label{eq:cc}
\varphi(A) = U^\ast A U \quad \mbox{for all} \quad A \in
B(\mathfrak{h})
\end{equation}
or
\begin{equation} \label{eq:dd}
\varphi(A) = U^\ast c^\ast A^\ast c U \quad \mbox{for all} \quad A
\in B(\mathfrak{h})
\end{equation}
where $c$ is the anti-unitary operator on $\mathfrak{h}$ induced
by complex conjugation. (As was shown in \cite[Theorem 5]{LMa} the
case where $\mathfrak{k} = \mathbb{C}$ and $\varphi$ is a vector
(pure) state also reduces to the form \ref{eq:cc} above. Just
select $x_\varphi \in \mathfrak{h}$ so that $\varphi = (\cdot
x_\varphi, x_\varphi)$ and define $U$ by $1 \rightarrow
x_\varphi$.)

If $\varphi$ is of the form \ref{eq:dd} above it may be regarded
as a map generated by the $*$-antiautomorphism $A \rightarrow c^*
A^* c$ on $B(\mathfrak{h})$ composed with the map $B(\mathfrak{h})
\rightarrow B(\mathfrak{k}): B \mapsto U^*BU$. Now by (\cite{W};
cf. \cite[Theorem 1]{Sh}) the dual of the $*$-antiautomorphism $A
\rightarrow c^* A^* c$ yields a bijection on the vector states of
$B(\mathfrak{h})$ which preserves transition probabilities. Thus
all that remains is to show that a map of the form \ref{eq:cc}
above preserves transition probabilities and has locally solid
range. Therefore assume that $\varphi$ is of this form. But then
$\varphi^*$ will map a vector state $$\omega(\cdot) = (\cdot
x_\omega, x_\omega) \quad x_\omega \in \mathfrak{k}$$ onto
$$\varphi^*(\omega)(\cdot) = (U*\cdot U x_\omega, x_\omega) =
(\cdot U x_\omega, Ux_\omega).$$ The injectivity of $U$ ensures
that $U^*U = \I$ and hence for any two unit vectors $x,y \in
\mathfrak{k}$ we have that $$|(x, y)|^2 = |(U^*Ux, y)|^2 = |(Ux,
Uy)|^2.$$ It therefore clearly follows that $\varphi^\sharp$
preserves transition probabilities.

Next let $E_\varphi = UU^*$ (the projection onto
$U(\mathfrak{k})$). Since $U^*$ restricts to a linear isometry
from $E_\varphi(\mathfrak{h}) = U(\mathfrak{k})$ onto
$\mathfrak{k}$, it is an exercise to show that $A \rightarrow
U^\ast A U$ and $A \rightarrow U^\ast c^\ast A^\ast c U$
respectively induce an onto $*$-isomorphism and an onto
$*$-anti-isomorphism from $B(E_\varphi(\mathfrak{h}))$ onto
$B(\mathfrak{k})$. Thus any map $\varphi$ of the form \ref{eq:cc}
above may be written as
\begin{equation} \label{eq:ee}
\varphi = \psi(E_\varphi\cdot E_\varphi)
\end{equation}
where $\psi$ is a $*$-isomorphism from
$B(E_\varphi(\mathfrak{h}))$ onto $B(\mathfrak{k})$. Now as we
noted in the discussion preceding this theorem, at the Hilbert
space level the set $\vee\{E_\rho| \rho = \varphi^\sharp(\omega),
\omega \mbox{ a vector state of }B(\mathfrak{k})\}$ corresponds to
nothing more than the closure of $\mathfrak{h}_0 =
\mathrm{span}\{\widetilde{x} \in \mathfrak{h} |
\omega_{\widetilde{x}} = \varphi^\sharp(\omega), \omega \mbox{ a
vector state of }B(\mathfrak{k}) \}$. However
$\mathrm{span}\{\widetilde{x} \in \mathfrak{h} |
\omega_{\widetilde{x}} = \varphi^\sharp(\omega), \omega \mbox{ a
vector state of }B(\mathfrak{k}) \}$ is of course precisely
$E_\varphi(\mathfrak{h}) = U(\mathfrak{k})$, and hence we may
conclude that here $$E_\varphi = \vee\{E_\rho| \rho =
\varphi^\sharp(\omega), \omega \mbox{ a vector state of
}B(\mathfrak{k})\}$$.

Now since $\psi$ is a bijective $*$-isomorphism we may apply
(\cite{W}; cf. \cite[Theorem 1]{Sh}) to conclude that the dual
$\psi^*$ maps the set of vector states of $B(\mathfrak{k})$ onto
that of $B(E_\varphi\mathfrak{h})$. Now any minimal projection
$E_0$ of $B(\mathfrak{h})$ with $E_0 \leq E_\varphi$ is of course
also a minimal projection of $B(E_\varphi\mathfrak{h})$. Hence
from what we have just shown, for such an $E_0$ we can always find
a minimal projection $F_0 \in B(\mathfrak{k})$ with
$$\psi^*(Tr_{\mathfrak{k}}(F_0\cdot F_0)) =
Tr_{E_\varphi\mathfrak{h}}(E_0\cdot E_0),$$ and hence by
\ref{eq:ee} above $$\varphi^*(Tr_{\mathfrak{k}}(F_0\cdot F_0)) =
Tr_{\mathfrak{h}}(E_0E_{\varphi}\cdot E_{\varphi}E_0) =
Tr_{\mathfrak{h}}(E_0\cdot E_0).$$ Thus $\varphi^\sharp$ has a
locally solid range.

\underline{$(2) \Rightarrow (1)$}: $\quad$ Let $\varphi^\sharp$ be
a transformation from the vector states of $B(\mathfrak{k})$ into
the vector states of $B(\mathfrak{k})$ which preserves transition
probabilities and has a locally solid range. By Lemma \ref{inj}
$\varphi^\sharp$ is injective. Now let $$E_\varphi = \vee\{E_\rho|
\rho = \varphi^\sharp(\omega), \omega \mbox{ a vector state of
}B(\mathfrak{k})\}.$$ We show that the dual of the map
$W_{E_\varphi}: B(\mathfrak{h}) \rightarrow
B(E_\varphi\mathfrak{h}): A \mapsto E_{\varphi}A E_{\varphi}$
bijectively maps the vector state space of
$B(E_\varphi\mathfrak{h})$ onto the range of $\varphi^\sharp$ in a
way that preserves transition probabilities. From this it then
follows that $\varphi^\sharp$ may be written as a bijection, say
$\psi^\sharp$, from the vector state space of $B(\mathfrak{k})$
onto that of $B(E_\varphi\mathfrak{h})$, composed with a
restriction of the dual of $W_{E_\varphi}$. The claim will then
follow from an application of (\cite{W}; cf. \cite[Theorem 1]{Sh})
to $\psi^\sharp$.

For any minimal projection $E_0$ of $B(\mathfrak{h})$ majorised by
$E_\varphi$, we have by hypothesis that $Tr_{\mathfrak{h}}(E_0
\cdot E_0)$ corresponds to an element of the range of
$\varphi^\sharp$. Moreover all elements of the range of
$\varphi^\sharp$ are of this form. These projections of course
correspond exactly to the minimal projections of
$B(E_\varphi\mathfrak{h})$. So for each such a minimal projection
$E_0 \in B(E_\varphi\mathfrak{h})$, the dual of $W_{E_\varphi}$
will map the vector state $Tr_{E_\varphi\mathfrak{h}}(E_0\cdot
E_0)$ of $B(E_\varphi\mathfrak{h})$ onto the corresponding element
$Tr_{\mathfrak{h}}(E_0\cdot E_0) \in B(\mathfrak{h})$ of the range
of $\varphi^\sharp$. Finally note that for any two minimal
projections $E,F \in B(\mathfrak{h})$ with $E, F \leq E_\varphi$,
we clearly have $$Tr_{E_\varphi\mathfrak{h}}(EF) =
Tr_{\mathfrak{h}}(E_\varphi EFE_\varphi) =
Tr_{\mathfrak{h}}(EF).$$ This fact can be shown to be equivalent
to the statement that $W_{E_\varphi}^*$ preserves transition
probabilities. The claim regarding $W_{E_\varphi}^*$ therefore
follows.

\underline{$(3) \Rightarrow (1)$}: $\quad$ The proof of this
implication is very similar to the proof of $(2) \Rightarrow (1)$
with the main difference being that here we use the previous
theorem instead of (\cite{W}; cf. \cite[Theorem 1]{Sh}). We
therefore leave this as an exercise.
\end{proof}

Based on the information garnered from the above result, a good
place to start searching for those transformations $\varphi^\sharp
: \mathcal{P}_{\mathcal{B}} \rightarrow \mathcal{P}_{\mathcal{A}}$
that are induced by linear maps on the underlying algebras, would
be among the co-orthogonal transformations that are locally
bi-orthogonal with a locally solid range. However for now we defer
such matters to the next section and content ourselves with the
following observation:

\begin{corollary}\label{cgwt}
Let $\mathcal{A, B}$ be $C^*$-algebras and $\varphi^\sharp$ a
transformation from $\mathcal{P}_{\mathcal{B}}$ into
$\mathcal{P}_{\mathcal{A}}$. If no irreducible representation of
$\mathcal{B}$ is of the form $M_2(\mathbb{C})$ and if
$\varphi^\sharp$ is fibre-preserving with locally solid range,
then $\varphi^\sharp$ is locally bi-orthogonal if and only if
$\varphi^\sharp$ locally preserves transition probabilities
\end{corollary}

\begin{proof} The ``if'' part is of course trivial and we therefore
indicate how the ``only if'' part may be verified. So suppose that
$\varphi^\sharp$ is fibre-preserving, locally bi-orthogonal, with
locally solid range. Let $\omega \in \mathcal{P}_{\mathcal{B}}$ be
given. By our supposition $\varphi^\sharp$ then restricts to a
bi-orthogonal transformation from $[\omega]$ into
$[\varphi^\sharp(\omega)]$. By Remark \ref{rar} the equivalence
classes $[\omega]$ and $[\varphi^\sharp(\omega)]$ respectively
correspond to the set of vector states of some $B(\mathfrak{k})$
and $B(\mathfrak{h})$, with $\mathcal{B}$ and $\mathcal{A}$
irreducibly represented on these two algebras. Call the induced
vector state map $\psi^\sharp$. Now by hypothesis $B(\mathfrak{k})
\neq M_2(\mathbb{C})$. So to be able to deduce the corollary from
the preceding theorem, all we need to do is to show that the local
solidness of the range of $\varphi^\sharp$ is enough to ensure the
local solidness of the range of $\psi^\sharp$ as defined in the
discussion preceding this theorem. This in turn is an exercise
depending on Remark \ref{rar} and Proposition \ref{orth}.
\end{proof}

\section{Pure state transformations induced by linear mappings}

We remind the reader that our primary challenge in this section is
to use the preceding results to identify those transformations
$\varphi^\sharp : \mathcal{P}_{\mathcal{B}} \rightarrow
\mathcal{P}_{\mathcal{A}}$ that are induced by linear maps on the
underlying algebras. We will do this in three phases. Suppose that
both $\mathcal{A}$ and $\mathcal{B}$ are reduced atomically
represented. Our first cycle of results will focus on identifying
those pure state transformations $\varphi^\sharp$ that correspond
to linear maps from $\mathcal{A}''$ into $\mathcal{B}''$. In the
second cycle we will show that under the assumption that
$\mathcal{A}$ and $\mathcal{B}$ are reduced atomically
represented, linear maps from $\mathcal{A}$ into $\mathcal{B}$
that have pure state preserving duals, live naturally in the class
of linear maps from $\mathcal{A}''$ into $\mathcal{B}''$. Hence
the correspondence established in the first cycle is therefore a
good place to start in searching for a solution to our primary
challenge. In the third and final cycle we will strengthen the
correspondence obtained in the first cycle by means of the
introduction of some additional continuity restrictions on
$\varphi^\sharp$ to finally obtain the solution to our problem.

Following a necessary technical lemma, we proceed to describe
those pure state transformations $\varphi^\sharp$ that correspond
to linear maps from $\mathcal{A}''$ into $\mathcal{B}''$.

\begin{lemma}\label{ps}
Suppose that $\mathcal{A}$ is a $C^*$-algebra which is reduced
atomically represented. Then the $\sigma$-weakly continuous pure
states of $\mathcal{A}''$ canonically correspond to the
pure states of $\mathcal{A}$.
\end{lemma}

\begin{proof}
The one direction of this correspondence follows from \cite[Lemma
4]{LMa} and the necessary normality ($\sigma$-weak continuity) of
the pure states of $\mathcal{A}$ in this representation (see
Remark \ref{rar}). For the converse recall that any normal
($\sigma$-weakly continuous) state $\omega$ of $\mathcal{A}''$ is
of the form $\omega = \sum_{k=1}^{\infty} \lambda_k \omega_{x_k}$
where $\lambda_k \geq 0$ for each $k$, $\sum_{k=1}^\infty
\lambda_k = 1$, and $\{x_k\} \subset \mathfrak{h} = \oplus_{a \in
\mathbb{A}} \mathfrak{h}_a$ is an orthonormal sequence
\cite[7.1.12]{KRi}. If therefore $\omega$ is to be an extreme
point of the state space of $\mathcal{A}''$, it may not be a
convex combination of distinct vector states and hence must itself
be a vector state, ie. $\lambda_k = 0$ for all but one $k$. Thus
suppose that $\omega = \omega_x$ for some norm-one $x \in
\oplus_{a \in \mathbb{A}} \mathfrak{h}_a$. Such an $x$ is of
course of the form $x = \oplus_{a \in \mathbb{A}} \mu_a x_a =
\oplus_{a \in \mathbb{A}} \mu_a \hat{x}_a$  with $\|\hat{x}_a\| =
1$ for each $a \in \mathbb{A}$ and $\sum_{a \in \mathbb{A}}
\abs{\mu_a}^2 = 1$. Now since each $\hat{x}_a$ is orthogonal to
all the elements of $\mathfrak{h}$ for which the $a^{\mathrm{th}}$
coordinate is zero, a careful consideration of the action of
$\omega = \omega_x$ on $\mathcal{A}'' = \oplus_{a \in \mathbb{A}}
B(\mathfrak{h}_a)$ reveals that it may be written in the form
$$\omega_x = \oplus_{a \in \mathbb{A}} \abs{\mu_a}^2 \omega_{x_a}
= \sum_{a \in \mathbb{A}} \abs{\mu_a}^2 \omega_{\hat{x}_a}.$$ As
before it is now clear that if $\omega = \omega_x$ is to be an
extreme point of the state space of $\mathcal{A}''$, we must have
that $\mu_a = 0$ for all but one $a \in \mathbb{A}$. From the
discussion in Remark \ref{rar} it is now clear that such an
$\omega$ canonically corresponds to a pure state of $\mathcal{A}$.
\end{proof}

\begin{theorem} \label{gen}
Let $\mathcal{A, B}$ be $C^\ast$-algebras with
$(\pi_{\mathcal{A}}, \mathfrak{h}_{\mathcal{A}})$ and
$(\pi_{\mathcal{B}}, \mathfrak{h}_{\mathcal{B}})$ their respective
reduced atomic representations, and let $\varphi^\sharp$ be a
transformation from $\mathcal{P}_{\mathcal{B}}$ into
$\mathcal{P}_{\mathcal{A}}$. Consider the following statements:
\begin{enumerate}
\item There exists a linear map $\widetilde{\varphi} :
\pi_{\mathcal{A}}(\mathcal{A})'' \rightarrow
\pi_{\mathcal{B}}(\mathcal{B})''$ such that
$\varphi^\sharp(\omega) = \omega \circ \widetilde{\varphi}$ for
every pure state $\omega \in \mathcal{P}_{\mathcal{B}}$.

\item $\varphi^\sharp$ is fibre-preserving, locally preserves
transition probabilities, and has a locally solid range.

\item $\varphi^\sharp$ is locally bi-orthogonal (and hence
fibre-preserving by Lemma \ref{coorth}) and has a locally solid
range.
\end{enumerate}
The implications $(1) \Leftrightarrow (2) \Rightarrow (3)$ hold in
general with all three statements being equivalent if no
irreducible representation of $\mathcal{B}$ is of the form
$M_2(\mathbb{C})$.
\end{theorem}

\begin{proof}
The implication $(2) \Rightarrow (3)$ is obvious in the light of
the results in section two. In addition the fact that $(3) \Rightarrow
(2)$ under the assumption that no irreducible representation of
$\mathcal{B}$ is of the form $M_2(\mathbb{C})$, is a direct consequence of
Lemma \ref{coorth} and Corollary \ref{cgwt}. Thus we need only prove that
$(1) \Leftrightarrow (2)$.

\textbf{\underline{$(1) \Rightarrow (2)$}:}$\quad$ Suppose that (1)
holds. For the sake of simplicity we may of course assume that both
$\mathcal{A}$ and $\mathcal{B}$ are reduced atomically represented. Then
$\mathcal{A}''$ and $\mathcal{B}''$ are of the form $\oplus_{a \in
\mathbb{A}} B(\mathfrak{h}_a)$ and $\oplus_{b \in \mathbb{B}}
B(\mathfrak{k}_b)$ respectively.

We first show that $\widetilde{\varphi}$ is necessarily a
contractive adjoint preserving map. To this end let $A = A^* \in
\oplus_{a \in \mathbb{A}} B(\mathfrak{h}_a)$ be given. A typical
$\sigma$-weakly continuous pure state $\omega$ of $\mathcal{B}''$ is
of course of the form $(\cdot\hat{z}_{b_0}, \hat{z}_{b_0})$ for some
unit vector $\hat{z}_{b_0} \in \oplus_{b \in \mathbb{B}} \mathfrak{k}_b$
(corresponding to some $z_{b_0} \in \mathfrak{k}_{b_0}$ -- see Remark
\ref{rar}). The transformation $\varphi^\sharp$ will then map this
pure state onto a similar looking $\sigma$-weakly continuous pure state
$\rho$ of $\mathcal{A}''$; say $(\cdot\hat{x}_{a_0}, \hat{x}_{a_0})$
where $x_{a_0} \in \mathfrak{h}_{a_0}$. Then by the assumption on
$\widetilde{\varphi}$ we have that $$(\widetilde{\varphi}(A)\hat{z}_{b_0},
\hat{z}_{b_0}) = (A\hat{x}_{a_0}, \hat{x}_{a_0}) \in \mathbb{R}$$ with
$$|(\widetilde{\varphi}(A)\hat{z}_{b_0}, \hat{z}_{b_0})| \leq \|A\|.$$
Since $b_0 \in \mathbb{B}$ and $z_{b_0} \in \mathfrak{k}_{b_0}$ was
arbitrary, this suffices to prove that $\widetilde{\varphi}$ preserves
adjoints with $\|\widetilde{\varphi}\| \leq 1$. In particular since
$\varphi^\sharp$ is induced by the dual of $\widetilde{\varphi}$, we
then have that $$\|\varphi^\sharp(\omega_0) - \varphi^\sharp(\omega_1)\|
\leq \|\omega_0 - \omega_1\|.$$ Thus $\varphi^\sharp$ must then be
co-orthogonal (since $\|\omega_0 - \omega_1\| < 2$ $\Rightarrow$
$\|\varphi^\sharp(\omega_0) - \varphi^\sharp(\omega_1)\| < 2$) and hence
fibre-preserving by Lemma \ref{coorth}.

It remains to show that $\varphi^\sharp$ locally preserves transition
probabilities and that it has a locally solid range. To this end let
$\omega_0 = (\cdot\hat{z}_{b_0}, \hat{z}_{b_0})$ be as before and
consider the action of $\widetilde{\varphi}^*$ on $[\omega_0]$. Recall
that $[\omega_0]$ corresponds to the vector states of $B(\mathfrak{k}_{b_0})$
in that the members of $[\omega_0]$ are precisely the vector states of
the form $(\cdot\hat{y}_{b_0}, \hat{y}_{b_0})$ for some
unit vector $y_{b_0} \in \mathfrak{k}_{b_0}$ -- see Remark
\ref{rar}). Now let $B \rightarrow E_0BE_0$ be the canonical compression
from $\mathcal{B}'' = \oplus_{b \in \mathbb{B}} B(\mathfrak{k}_b)$ onto
$B(\mathfrak{k}_{b_0})$. It is then an exercise to show that the action
of $\widetilde{\varphi}^*$ from $[\omega_0]$ into
$\mathcal{P}_{\mathcal A}$ corresponds canonically to the action induced
by $E_0\widetilde{\varphi}E_0$ on the vector states of
$B(\mathfrak{k}_{b_0})$. If necessary we may therefore replace
$\widetilde{\varphi}$ by $E_0\widetilde{\varphi}E_0$, and assume that
$\mathcal{B}'' = B(\mathfrak{k})$ where now the vector states of
$B(\mathfrak{k})$ plays the role of $[\omega_0]$.

Now since $\widetilde{\varphi}^*$ is fibre-preserving, all the
vector states of $B(\mathfrak{k})$ get mapped onto a single
equivalence class $[\varphi^\sharp(\omega_0)]$ of $\mathcal{A}'' =
\oplus_{a \in \mathbb{A}} B(\mathfrak{h}_a)$. This single
equivalence class corresponds to the vector states of precisely
one of the $B(\mathfrak{h}_a)$'s; say $B(\mathfrak{h}_{a_0})$. So
in the notation of Remark \ref{rar}, this means that for every
vector state $\omega_z = (\cdot z, z)$ of $B(\mathfrak{k})$ we
have $$\varphi^\sharp(\omega_z) = (\cdot\hat{x}_{a_0},
\hat{x}_{a_0})$$ for some unit vector $x_{a_o} \in
\mathfrak{h}_{a_0}$. In particular this means that for any unit
vector $z \in B(\mathfrak{k})$ and any $A = \oplus_{a \in
\mathbb{A}} A_a \in \mathcal{A}'' = \oplus_{a \in \mathbb{A}}
B(\mathfrak{h}_a)$ we have
\begin{equation} \label{eq:az}
(\widetilde{\varphi}(A)z, z) = (A\hat{x}_{a_0},
\hat{x}_{a_0}) = (A_{a_0}x_{a_0}, x_{a_0}).
\end{equation}
Clearly $\widetilde{\varphi}$ annihilates all elements $A =
\oplus_{a \in \mathbb{A}} A_a$ of $\mathcal{A}''$ with $0$ in the
${a_0}^{\mathrm{th}}$ coordinate. More to the point if by
$\hat{A}_{a_0}$ we denote the element of $\mathcal{A}''$ with
$A_{a_0}$ in the ${a_0}^{\mathrm{th}}$ coordinate and zeros
elsewhere, this shows that $\widetilde{\varphi}$ factors through
$B(\mathfrak{h}_{a_0})$ in that we may write it as a composition
of the maps  $$\mathcal{A}'' = \oplus_{a \in \mathbb{A}}
B(\mathfrak{h}_a) \rightarrow B(\mathfrak{h}_{a_0}): \oplus_{a \in
\mathbb{A}} A_a \mapsto A_{a_0}$$ and $$\widetilde{\varphi}_0 :
B(\mathfrak{h}_{a_0}) \rightarrow B(\mathfrak{k}): A_{a_0}
\rightarrow \widetilde{\varphi}(\hat{A}_{a_0}).$$ It is now clear
from equation (\ref{eq:az}) that $\widetilde{\varphi}_0^*$ maps
the vector states of $B(\mathfrak{k})$ into the vector states of
$B(\mathfrak{h}_{a_0})$ in a way that canonically agrees with the
action of $\widetilde{\varphi}^*$ from $[\omega_0]$ to
$[\varphi^\sharp(\omega_0)]$. We may therefore apply Theorem
\ref{gwt} to see that $\widetilde{\varphi}_0^*$ preserves
transition probabilities and that (by abuse of notation for the
sake of clarity)
$$\widetilde{\varphi}_0^*([\omega_0])^{\perp\perp} =
\widetilde{\varphi}_0^*([\omega_0])$$ (complements taken with
respect to $B(\mathfrak{h}_{a_0})$). It is then not difficult to
see that this is the same as saying that $\widetilde{\varphi}^*$
preserves transition probabilities on $[\omega_0]$ and that
$$\widetilde{\varphi}^*([\omega_0])^{\perp\perp} =
\widetilde{\varphi}^*([\omega_0])$$ (complements taken with
respect to $\mathcal{A}'' = \oplus_{a \in \mathbb{A}}
B(\mathfrak{h}_a)$). In the light of Proposition \ref{orth} and
the discussion preceding Theorem \ref{gwt}, this last statement is
just another way of saying that for a linear subspace S of
$\mathfrak{h}_{a_0}$ the claim $S^{\perp\perp} = S$ (comlements
taken in $\mathfrak{h}_{a_0}$) is equivalent to the claim
$\{\hat{x}_{a_0} | x_{a_0} \in S\}^{\perp\perp} = \{\hat{x}_{a_0}
| x_{a_0} \in S\}$ (complements taken in $\oplus_{a \in
\mathbb{A}} \mathfrak{h}_a$. Since our original choice of
$\omega_0$ was arbitrary, this proves the required implication.

\textbf{\underline{$(2) \Rightarrow (1)$}:}$\quad$ Suppose that
(2) holds. By Remark \ref{rar} we have that $$\mathcal{A}'' =
\oplus_{a \in \mathbb{A}} B(\mathfrak{h}_a) \quad \mbox{and} \quad
\mathbb{B}'' = \oplus_{b \in \mathbb{B}} B(\mathfrak{k}_b)$$ with
the classes of vector states of each distinct $B(\mathfrak{h}_a)$
and $B(\mathfrak{k}_b)$ corresponding in a unique way to distinct
equivalence classes of pure states of $\mathcal{A}$ and
$\mathcal{B}$ respectively. Now since $\varphi^\sharp$ is
fibre-preserving, we may re-index the $B(\mathfrak{k}_b)$'s with a
double index where for any $a \in \mathbb{A}$, the collection
$B(\mathfrak{k}_{\lambda}^{(a)})$ $(\lambda \in \Lambda_a)$
denotes all the $B(\mathfrak{k}_b)$'s that correspond to
equivalence classes of pure states of $\mathcal{B}$ that map into
the single equivalence class of $\mathcal{P}_{\mathcal A}$
corresponding to $B(\mathfrak{h}_a)$. (Note that some of the
$\Lambda_a$'s may be empty.) It then surely follows that
$$\mathcal{B}'' = \oplus_{a \in \mathbb{A}} (\oplus_{\lambda \in
\Lambda_a} B(\mathfrak{k}_{\lambda}^{(a)}).$$ Now fix a non-empty
$\Lambda_a$ and let $\lambda \in \Lambda_a$ be given. Since the
action of $\varphi^\sharp$ from the equivalence class
corresponding to $B(\mathfrak{k}_{\lambda}^{(a)})$ into the one
corresponding to $B(\mathfrak{h}_a)$ preserves transition
probabilities and has a locally solid range, we may apply Theorem
\ref{gwt} to obtain a contractive linear map
$\widetilde{\varphi}^{(a)}_\lambda : B(\mathfrak{h}_a) \rightarrow
B(\mathfrak{k}_{\lambda}^{(a)})$ whose dual in a canonical way
induces the action of $\varphi^\sharp$ on the corresponding
equivalence classes of pure states. Summing over $\lambda$ we get
a contractive map $$\widetilde{\varphi}_a : B(\mathfrak{h}_a)
\rightarrow \oplus_{\lambda \in \Lambda_a}
B(\mathfrak{k}_{\lambda}^{(a)}) : A_a \mapsto \oplus_{\lambda \in
\Lambda_a} \widetilde{\varphi}^{(a)}_\lambda(A_a)$$ whose dual
again in a canonical way induces the action of $\varphi^\sharp$
from all the equivalence classes corresponding to the
$B(\mathfrak{k}_{\lambda}^{(a)})$'s $(\lambda \in \Lambda_a)$ to
the one corresponding to $B(\mathfrak{h}_a)$. Summing over
$\mathbb{A}_0 = \{a \in \mathbb{A} | \Lambda_a \neq \emptyset\}$,
we get a contractive map $$\widetilde{\varphi} = \oplus_{a \in
\mathbb{A}_0} \widetilde{\varphi}_a : \oplus_{a \in \mathbb{A}_0}
B(\mathfrak{h}_a) \rightarrow \oplus_{a \in
\mathbb{A}_0}(\oplus_{\lambda \in \Lambda_a}
B(\mathfrak{k}_{\lambda}^{(a)})) =\mathcal{B}''$$ which may be
extended to all of $\mathcal{A}'' = \oplus_{a \in \mathbb{A}}
B(\mathfrak{h}_a)$ by defining it in such a way that it
annihilates all the $B(\mathfrak{h}_a)$'s for which the
corresponding $\Lambda_a$ is empty. By construction the dual of
this extended map induces the action of $\varphi^\sharp$. Thus (1)
holds.
\end{proof}

We list the details of the bijective case separately because of
its more elegant behaviour. For this case we note that
while the statement regarding bi-orthogonal bijections is new, the
case pertaining to bijections which preserve transition
probabilities is basically just a version of Shultz's result (see
Theorem 1.5) with the condition regarding orientation removed.

\begin{proposition}\label{genbij}
Let $\mathcal{A, B}$ be $C^\ast$-algebras with
$(\pi_{\mathcal{A}}, \mathfrak{h}_{\mathcal{A}})$ and
$(\pi_{\mathcal{B}}, \mathfrak{h}_{\mathcal{B}})$ their respective
reduced atomic representations.

For any Jordan $\ast$-isomorphism $\widetilde{\varphi}$ from
$\pi_{\mathcal{A}}(\mathcal{A})''$ onto
$\pi_{\mathcal{B}}(\mathcal{B})''$ the adjoint of
$\widetilde{\varphi}$ restricts to bijection
$\varphi^\sharp$ from the set of pure states of $\mathcal{B}$ onto
the set of pure states of $\mathcal{A}$ which preserves all the
transition probabilities . In particular $\varphi^\sharp$ is
bi-orthogonal.

Conversely if either $\varphi^\sharp$ is a bijection from the pure
states of $\mathcal{B}$ onto the pure states of $\mathcal{A}$
which preserves transition probabilities or if
$\varphi^\sharp$ is a bi-orthogonal bijection and either $\mathcal{A}$
or $\mathcal{B}$ has the property that no irreducible representation
is of the form $M_2(\mathbb{C})$, then we can find a Jordan
$\ast$-isomorphism $\widetilde{\varphi}$ from
$\pi_{\mathcal{A}}(\mathcal{A})''$ onto
$\pi_{\mathcal{B}}(\mathcal{B})''$ such that
$$\varphi^\sharp(\omega)(\cdot) =
\omega(\pi_{\mathcal{B}}^{-1}\widetilde{\varphi}(\cdot)\pi_{\mathcal{A}})$$
for each pure state $\omega$ of $\mathcal{B}$.

(In the above statements we have again identified the pure states
of $\mathcal{A}$ and $\mathcal{B}$ with their unique normal
extensions to $\pi_{\mathcal{A}}(\mathcal{A})''$ and
$\pi_{\mathcal{B}}(\mathcal{B})''$ \cite[Lemma 11]{LMa}.)
\end{proposition}

\begin{proof}
In the proof we concentrate on the case pertaining to
bi-orthogonality. The proofs of the two cases are similar with
Wigner's result (\cite{W}; cf. \cite[Theorem 1]{Sh}) being used in
the former case instead of Theorem \ref{bijbh} (\cite{CdVLL}).

Without any loss of generality we may identify both $\mathcal{A}$
and $\mathcal{B}$ with their respective reduced atomic
representations.

First let $\varphi^\sharp$ be an orthogonal bijection from the set
of pure states of $\mathcal{B}$ onto the set of pure states of
$\mathcal{A}$ and suppose that no irreducible representation of
$\mathcal{A}$ is of the form $M_2(\mathbb{C})$. Now by assumption
$\mathcal{A}''$ and $\mathcal{B}''$ are of the form $\oplus_{a \in
\mathbb{A}} B(\mathfrak{h}_a)$ and $\oplus_{b \in \mathbb{B}}
B(\mathfrak{k}_b)$ respectively, with each equivalence class of
pure states of $\mathcal{A}$ and $\mathcal{B}$ corresponding to
the set of vector states of one of the $B(\mathfrak{h}_a)$'s and
$B(\mathfrak{k}_b)$'s respectively in the sense described in
Remark \ref{rar}. Now by Lemma \ref{coorth} both $\varphi^\sharp$
and its inverse is fibre-preserving. It therefore induces a
bijection from the set of equivalence classes of pure states of
$\mathcal{B}$ onto the set of equivalence classes of pure states
of $\mathcal{A}$. We may therefore re-index the $B(\mathfrak{k}_b)$'s
($b \in \mathbb{B}$) with the index set $\mathbb{A}$ in such a way
that $\mathcal{B}'' = \oplus_{a \in \mathbb{A}} B(\mathfrak{k}_a)$,
with $\varphi^\sharp$ for each $a \in \mathbb{A}$ mapping the equivalence
class corresponding to the vector states of $B(\mathfrak{k}_a)$ onto the
equivalence class corresponding to the vector states of $B(\mathfrak{h}_a)$.

Now by the assumption on $\mathcal{A}$ we have that
$B(\mathfrak{h}_a) \neq M_2(\mathbb{C})$ for each $a \in
\mathbb{A}$. Since $\varphi^\sharp$ is also bi-orthogonal, it now
follows from Proposition \ref{orth} and Theorem \ref{bijbh}
(\cite{CdVLL}) that for each $a \in \mathbb{A}$, $\varphi^\sharp$
induces either a $\ast$-isomorphism or $\ast$-antiisomorphism
$\widetilde{\varphi}_a$ from $B(\mathfrak{h}_a)$ onto
$B(\mathfrak{k}_a)$. Moreover for each $a \in \mathbb{A}$ the
transformation that $\varphi^\sharp$ induces from the vector
states of $B(\mathfrak{k}_a)$ onto the vector states of
$B(\mathfrak{h}_a)$ appears as a restriction of the adjoint of
$\widetilde{\varphi}_a$. Recall that in the notation of Remark
\ref{rar} pure states of $\mathcal{B}$ correspond to vector states
of the form $\omega_{\hat{x}_a}$ for some $a \in \mathbb{A}$ and
some norm one vector $x_a \in \mathfrak{k}_a$. With this in mind
it is now a simple matter to verify that $\widetilde{\varphi} =
\oplus_{a \in \mathbb{A}} \widetilde{\varphi}_a$ is a Jordan
$\ast$-isomorphism from $\mathcal{A}'' = \oplus_{a \in \mathbb{A}}
B(\mathfrak{h}_a)$ onto $\mathcal{B}'' = \oplus_{a \in \mathbb{A}}
B(\mathfrak{k}_a)$ with the property that
$$\varphi^\sharp(\omega_{\hat{x}_a}) = \omega_{\hat{x}_a} \circ
\widetilde{\varphi}$$ for each $a \in \mathbb{A}$ and each norm
one vector $x_a \in \mathfrak{k}_a$. This clearly suffices to
establish the claim.

Conversely let $\widetilde{\varphi}$ be a Jordan
$\ast$-isomorphism from $\mathcal{A}''$ onto $\mathcal{B}''$. By
Lemma 2.2 $\widetilde{\varphi}$ is then necessarily a
$\sigma$-weak homeomorphism.  Moreover by \cite[Theorem 5]{LMa}
the adjoint of $\widetilde{\varphi}$ restricts to a bijection from
the pure states of $\mathcal{B}''$ onto the pure states of
$\mathcal{A}''$. Thus $\widetilde{\varphi}^\ast$ trivially maps
the $\sigma$-weakly continuous pure states of $\mathcal{B}''$ onto
the $\sigma$-weakly continuous pure states of $\mathcal{A}''$.

It remains to show that the restriction of
$\widetilde{\varphi}^\ast$ to the pure states of $\mathcal{B}$
preserves transition probabilities. We once again identify the
pure states of $\mathcal{A}$ and $\mathcal{B}$ with their unique
$\sigma$-weakly continuous extensions to $\mathcal{A}''$ and
$\mathcal{B}''$. The fact that $\varphi^\sharp$ is bi-orthogonal
is a straightforward consequence of the fact that
$\widetilde{\varphi}$, and hence also its dual, is a surjective
linear isometry \cite[3.2.3]{BRo}. Thus given any two pure states
$\omega_0$ and $\omega_1$ of $\mathcal{B}$ it easily follows that
$\|\omega_0 - \omega_1\| = 2$ if and only if
$\|\varphi^\sharp(\omega_0) - \varphi^\sharp(\omega_1)\| =
\|\omega_0 \circ \widetilde{\varphi} - \omega_1 \circ
\widetilde{\varphi}\| = 2$.

Establishing that $\varphi^\sharp$ preserves transition
probabilities requires a bit more work and may in fact be deduced
from Wigner's classical result. We show how to do this using the
available structure. As before bi-orthogonality ensures that both
$\varphi^\sharp$ and its inverse are fibre-preserving. Hence it is
enough to prove the preservation of transition probabilities for
sets of unitarily equivalent pure states only.  Recall that
$\mathcal{A}''$ and $\mathcal{B}''$ are of the form $\mathcal{A}''
= \oplus_{a \in \mathbb{A}} B(\mathfrak{h}_a)$ and $\mathcal{B}''
= \oplus_{b \in \mathbb{B}} B(\mathfrak{k}_b)$. The fact that both
$\varphi^\sharp$ and its inverse are fibre-preserving ensures that
we may re-index the expression for $\mathcal{B}''$ with the index
set $\mathbb{A}$ in such a way that for any given $a \in
\mathbb{A}$, $\varphi^\sharp$ identifies the pure states
corresponding to the vector states of $B(\mathfrak{k}_a)$ with the
pure states corresponding to the vector states of
$B(\mathfrak{h}_a)$ (see Remark \ref{rar}). For any $B = \oplus_{a
\in \mathbb{A}} B_a \in \oplus_{a \in \mathbb{A}}
B(\mathfrak{h}_a)$ and any given $d \in \mathbb{A}$, we then have
that $\omega_{\hat{x}_d}(B) = 0$ for all unit vectors $x_d \in
\mathfrak{h}_d$ if and only if
$\omega_{\hat{y}_d}(\widetilde{\varphi}(B)) = 0$ for all unit
vectors $y_d \in \mathfrak{k}_d$ (again notation is as in Remark
\ref{rar}). It is now an exercise to see that this ensures that
for each $a \in \mathbb{A}$, $\widetilde{\varphi}$ maps the copy
of $B(\mathfrak{h}_a)$ in $\mathcal{A}''$ onto the copy of
$B(\mathfrak{k}_a)$ in $\mathcal{B}''$. All that remains to be
done is to check that the map that $\varphi^\sharp$ induces from
the vector states of $B(\mathfrak{k}_a)$ onto those of
$B(\mathfrak{h}_a)$, is a restriction of the dual of the map
$\widetilde{\varphi}$ induces from $B(\mathfrak{h}_a)$ to
$B(\mathfrak{k}_a)$. Since the map induced by
$\widetilde{\varphi}$ is necessarily either a $*$-isomorphism or a
$*$-anti-isomorphism \cite[3.2.2]{BRo}, we may then directly apply
Wigner's theorem (\cite{W}; cf, \cite[Theorem 1]{Sh}) to get the
required conclusion.
\end{proof}

The above results provide the context for a non-commutative
version of the Banach-Stone theorem \cite[3.4.3]{KRi}. To see
this recall that for any compact Hausdorff set $K$ the pure state
space of $C(K)$ corresponds exactly to the set of point
evaluations engendered by elements of $K$. Endowed with the
$\mathrm{weak}^\ast$ topology, this set of point evaluations is of
course homeomorphic to $K$ itself. The link between pure state
transformations and linear maps described above is therefore in the
same spirit as the link between continuous transformations on
compact Hausdorff sets and $\ast$-homomorphisms on $C(K)$-spaces.

In further support of this contention we present the next step in
our programme, which is to identify the linear maps from $\mathcal{A}$
into $\mathcal{B}$ whose duals preserve pure states, as a subclass of the
$\sigma$-weakly continuous linear maps from $\mathcal{A}''$ into
$\mathcal{B}''$ whose duals preserve $\sigma$-weakly continuous pure
states.

\begin{proposition}\label{extn}
As before let $\mathcal{A, B}$ be $C^\ast$-algebras with
$(\pi_{\mathcal{A}}, \mathfrak{h}_{\mathcal{A}})$ and
$(\pi_{\mathcal{B}}, \mathfrak{h}_{\mathcal{B}})$ their respective
reduced atomic representations. For any linear map $\varphi$ from
$\mathcal{A}$ into $\mathcal{B}$ with a pure state preserving dual,
$\pi_{\mathcal{B}} \circ \varphi \circ {\pi_{\mathcal{A}}}^{-1}$
is $\sigma$-weakly continuous and admits of a unique extension to
a $\sigma$-weakly continuous linear map $\widetilde{\varphi}$ from
$\pi_{\mathcal{A}}(\mathcal{A})''$ into
$\pi_{\mathcal{B}}(\mathcal{B})''$ whose dual preserves $\sigma$-weakly
continuous pure states.
In particular if $\varphi$ is a Jordan $*$-isomorphism from $\mathcal{A}$
onto $\mathcal{B}$, then its dual restricts to a bijection between the
respective sets of pure states, and the canonical
extension $\widetilde{\varphi}$ described above is a Jordan $*$-isomorphism
from $\pi_{\mathcal{A}}(\mathcal{A})''$ onto
$\pi_{\mathcal{B}}(\mathcal{B})''$.
\end{proposition}

\begin{proof}
Without loss of generality let $\mathcal{A, B}$ be identified with
their reduced atomic representations and let $\varphi: \mathcal{A}
\rightarrow \mathcal{B}$ be a linear map with a pure state
preserving dual. It follows from \cite[Theorem 5]{LMa} that
$\varphi^\ast$ is necessarily contractive. The map $\varphi^\ast$
is then bounded and so restricts to an affine map from the
norm-closed convex hull of the pure states of $\mathcal{B}$ into
the norm-closed convex hull of the pure states of $\mathcal{A}$.
We show that these norm closed convex hulls are precisely the
normal state spaces of $\mathcal{B}$ and $\mathcal{A}$
respectively. By Remark \ref{rar} all the pure states of say
$\mathcal{A}$ are necessarily normal and hence the norm-closed
convex hull of the pure states is at least contained in the normal
state space. (The fact that the normal state space of a concrete
$C^\ast$-algebra is a norm-closed convex set follows from for
example \cite[10.1.15]{KRi}.) Since by \cite[7.1.12 \&
10.1.11(i)]{KRi} the normal state space of say $\mathcal{A}$ is
the norm-closed convex hull of the vector states, we need only
show that the vector states are contained in the norm-closed
convex hull of the pure states to conclude that the two sets are
equal. This in turn can be seen as follows: Given a norm-one
element $x \in \mathfrak{h}_{\mathcal{A}}$, we saw in the proof of
Lemma \ref{ps} that in this case $\omega_x$ can be written in the
form $$\omega_x = \sum_{a \in \mathbb{A}} \abs{\mu_a}^2
\omega_{\hat{x}_a}$$ where each $\omega_{\hat{x}_a}$ is a pure
state and $\sum_{a \in \mathbb{A}} \abs{\mu_a}^2 = 1$. It follows
that $\varphi^\ast$ affinely maps the normal state space of
$\mathcal{B}$ into the normal state space of $\mathcal{A}$.

Now since the sets of $\sigma$-weakly continuous states of
$\mathcal{B}$ and $\mathcal{A}$ correspond canonically to the sets
of $\sigma$-weakly continuous states (ie. the normal states) of
$\mathcal{B}''$ and $\mathcal{A}''$ respectively (see
\cite[10.1.11(i)]{KRi}), it follows that the linear span of the
sets of $\sigma$-weakly continuous states of $\mathcal{A}$ and
$\mathcal{B}$ correspond canonically to $(\mathcal{A}'')_*$ and
$(\mathcal{B}'')_*$ respectively. Let $\widetilde{\varphi}$ be the
($\sigma$-weakly continuous) dual of the map induced by
$\varphi^*$ from $(\mathcal{B}'')_*$ to $(\mathcal{A}'')_*$. Then
for any normal state $\omega$ of $\mathcal{B}''$ and any $A \in
\mathcal{A}$ we have by construction that
$$\omega(\widetilde{\varphi}(A)) = {\widetilde{\varphi}}_\ast
\circ \omega(A) = \varphi^\ast \circ \omega(A) =
\omega(\varphi(A)).$$ Thus $\widetilde{\varphi}$ is
$\sigma$-weakly continuous, and since we must have that
$\widetilde{\varphi}(A) = \varphi(A)$ for each $A$, it is of
course also an extension of $\varphi$. The claim regarding the
uniqueness of the extension follows from for example
\cite[10.1.10]{KRi}. The claim regarding the $\sigma$-weakly
continuous pure states is then a trivial consequence of what we
just proved, and Lemma \ref{ps}.

Now if $\varphi$ is a Jordan $*$-isomorphism from $\mathcal{A}$
onto $\mathcal{B}$, then by \cite[Theorem 5]{LMa} its dual defines
a bijection between the respective sets of pure states. We may
then mimic the above argument to in this case construct an affine
bijection between the normal state spaces of $\mathcal{B}''$ and
$\mathcal{A}''$. By Kadison's result (cf. Theorem 1.3(1)) this
affine bijection yields a Jordan $*$-isomorphism
$\widetilde{\varphi}$ from $\mathcal{A}''$ onto $\mathcal{B}''$
which by the same argument as before can be shown to be the
required extension.
\end{proof}

We are now finally ready to present our main theorem. In
mathematical terms this amounts to a very general non-commutative
Banach-Stone type theorem. At a slightly different level one may
interpret the bijective case of this as a result stating that
there is indeed enough information internally encoded in the pure
state spaces of two $C^\ast$-algebras $\mathcal{A}$ and
$\mathcal{B}$ to ensure their physical equivalence via a suitable
pure state bijection and also enough to be able to identify those
pure state bijections which actually correspond to some type of
Wigner symmetry. We also take this opportunity to invite the
reader to compare the result below with Theorem 5.7 of \cite{Sto}.
In this result St{\o}rmer uses conditions analogous to
\emph{fibre-preserving} and \emph{locally solid} to describe a
class of Jordan $*$-homomorphisms in terms of their action on pure
states.

We point out that we have assumed our algebras to be unital. For the
result to hold in the non-unital case $\varphi^\sharp$ is required to
admit of a homeomorphic action from $\mathcal{P}_{\mathcal B} \cup
\{0\}$ to $\mathcal{P}_{\mathcal A} \cup \{0\}$ which fixes 0 (see
\cite{Sh} and \cite{Br}). We take this opportunity
to re-emphasise the fact noted prior to Theorem \ref{bijbh} that the
restriction regarding $M_2(\mathbb{C})$ can not be removed.

\begin{theorem}[Non-commutative Banach- Stone theorem]\label{ncbs}
Let $\mathcal{A, B}$ be $C^\ast$-algebras and let $\varphi^\sharp$ be
a transformation from $\mathcal{P}_{\mathcal{B}}$ into
$\mathcal{P}_{\mathcal{A}}$. Consider the following
statements:
\begin{enumerate}
\item There exists a linear map $\varphi : \mathcal{A} \rightarrow
\mathcal{B}$ such that $\varphi^\sharp(\omega) = \omega \circ
\varphi$ for every pure state $\omega \in \mathcal{P}_{\mathcal{B}}$.
\item $\varphi^\sharp$ is uniformly $\sigma(\mathcal{B}^\ast, \mathcal{B}) -
\sigma(\mathcal{A}^\ast, \mathcal{A})$ continuous, fibre-preserving,
locally preserves transition probabilities, and has a locally solid range.
\item $\varphi^\sharp$ is uniformly $\sigma(\mathcal{B}^\ast, \mathcal{B}) -
\sigma(\mathcal{A}^\ast, \mathcal{A})$ continuous, locally bi-orthogonal (and
hence fibre-preserving by Lemma \ref{coorth}) and has a locally solid
range.
\end{enumerate}
The implications $(1) \Leftrightarrow (2) \Rightarrow (3)$ hold in
general with all three statements being equivalent if no
irreducible representation of $\mathcal{B}$ is of the form
$M_2(\mathbb{C})$.
\end{theorem}

\begin{proof} The implication $(2) \Rightarrow (3)$ is fairly clear, whereas
the fact that $(3) \Rightarrow (2)$ whenever no
irreducible representation of $\mathcal{B}$ is of the form
$M_2(\mathbb{C})$, follows from Corollary \ref{cgwt}. Hence we need only show
that $(1) \Leftrightarrow (2)$.

Firstly let $\varphi$ be a linear map from $\mathcal{A}$ into
$\mathcal{B}$ with a pure state preserving dual. Then all
statements in (2) except the claim about uniform weak* continuity
follows from a combination of Theorem \ref{gen} and Proposition
\ref{extn}. To see the last claim note that since $\varphi$ is
necessarily continuous (see \cite[Theorem 5]{LMa}), its dual is of
course $\sigma(\mathcal{B}^\ast, \mathcal{B}) -
\sigma(\mathcal{A}^\ast, \mathcal{A})$ continuous on the
$\sigma(\mathcal{B}^\ast, \mathcal{B})$--compact unit ball of
$\mathcal{B}^*$. Continuity of the dual on a compact superset then
ensures that the restriction to the subset
$\mathcal{P}_\mathcal{B}$ must necessarily be uniformly
$\sigma(\mathcal{B}^\ast, \mathcal{B}) - \sigma(\mathcal{A}^\ast,
\mathcal{A})$ continuous.

For the converse we may assume without loss of generality that
both $\mathcal{A}$ and $\mathcal{B}$ are reduced atomically
represented. It is then clear from Theorem \ref{gen} that the
hypotheses are sufficient to guarantee that $\varphi^\sharp$ is
induced by a linear map $\widetilde{\varphi}$ from $\mathcal{A}''$
into $\mathcal{B}''$. We therefore need only show that requiring
$\varphi^\sharp$ to in addition be uniformly $\textrm{weak}^\ast$
continuous is sufficient to guarantee that
$\widetilde{\varphi}(\mathcal{A}) \subset \mathcal{B}$.

Hence suppose that $\widetilde{\varphi}^\ast$ is uniformly
$\sigma(\mathcal{B}^\ast, \mathcal{B}) - \sigma(\mathcal{A}^\ast,
\mathcal{A})$ continuous from $\mathcal{P}_{\mathcal B}$ into
$\mathcal{P}_{\mathcal A}$. Recall the the biduals of
$\mathcal{A}$ and $\mathcal{B}$ may be identified with the double
commutants of their respective universal representations
\cite[10.1.1]{KRi}. Since in the universal representation of a
$C^\ast$-algebra all states are normal, it follows from
\cite[Proposition 6]{LMa} that each of $\mathcal{A}^{\ast\ast}$
and $\mathcal{B}^{\ast\ast}$ admit of central projections
$E_{\mathcal A}$ and $E_{\mathcal B}$ respectively such that
$\mathcal{A}_{E_{\mathcal A}}$ and $\mathcal{B}_{E_{\mathcal B}}$
are $\sigma$-weakly $\ast$-isomorphic to the reduced atomic
representations of $\mathcal{A}$ and $\mathcal{B}$ respectively.
In fact these projections correspond to nothing more than the
canonical central projections $z_{\mathcal{A}}$ and
$z_{\mathcal{B}}$ of $\mathcal{A}^{**}$ and $\mathcal{B}^{**}$
onto their respective atomic parts (as defined and used by both
Shultz \cite{Sh} and Brown \cite{Br}). Again in the notation of
Brown, $\widetilde{\varphi}$ must then correspond to a linear map
from $z_{\mathcal{A}}{\mathcal A}^{\ast\ast}$ to
$z_{\mathcal{B}}{\mathcal B}^{\ast\ast}$. Therefore given any $A
\in \mathcal{A}$, $\widetilde{\varphi}(A)$ corresponds to an
element $B$ of $z{\mathcal B}^{\ast\ast}$. Now since $A$ defines a
linear $\sigma(\mathcal{A}^\ast, \mathcal{A})$ continuous
functional on $\mathcal{A}^\ast$ and $\widetilde{\varphi}^*$ is
uniformly $\sigma(\mathcal{B}^\ast, \mathcal{B}) -
\sigma(\mathcal{A}^\ast, \mathcal{A})$ continuous from
$\mathcal{P}_{\mathcal B}$ into $\mathcal{P}_{\mathcal A}$, it
follows that $\widetilde{\varphi}(A)$ (that is $B$) is uniformly
$\sigma(\mathcal{B}^\ast, \mathcal{B})$ continuous on
$\mathcal{P}_{\mathcal B}$. Thus by \cite[Corollary 8]{Br} we must
have that $B \in z{\mathcal B}$. Given that $z{\mathcal B}$ (that
is $\mathcal{B}_{E_{\mathcal B}}$) corresponds to the reduced
atomic representation of $\mathcal{B}$ and that $B$ is identified
with $\widetilde{\varphi}(A)$ under this correspondence, this
shows that $\widetilde{\varphi}(A) \in \mathcal{B}$ as required.
\end{proof}

\begin{corollary}\label{ncbsbij}
Let $\mathcal{A, B}$ be $C^\ast$-algebras and let $\varphi^\sharp$ be
a transformation from $\mathcal{P}_{\mathcal{B}}$ into
$\mathcal{P}_{\mathcal{A}}$. Consider the following
statements:
\begin{enumerate}
\item There exists a Jordan $*$-isomorphism $\varphi$ from $\mathcal{A}$ onto
$\mathcal{B}$ such that $\varphi^\sharp(\omega) = \omega \circ
\varphi$ for every pure state $\omega \in \mathcal{P}_{\mathcal{B}}$.
\item $\varphi^\sharp$ is a bijective uniform $\sigma(\mathcal{B}^\ast,
\mathcal{B}) - \sigma(\mathcal{A}^\ast, \mathcal{A})$ homeomorphism which
preserves transition probabilities.
\item $\varphi^\sharp$ is bijective and is a biorthogonal uniform
$\sigma(\mathcal{B}^\ast, \mathcal{B}) - \sigma(\mathcal{A}^\ast,
\mathcal{A})$ homeomorphism.
\end{enumerate}
The implications $(1) \Leftrightarrow (2) \Rightarrow (3)$ hold in
general with all three statements being equivalent if either $\mathcal{A}$
or $\mathcal{B}$ has the property that no irreducible representation is of
the form $M_2(\mathbb{C})$.
\end{corollary}

\begin{proof}
First recall that any Jordan $*$-isomorphism $\varphi$ from
$\mathcal{A}$ onto $\mathcal{B}$ is a linear isometry
\cite[3.2.3]{BRo}. Hence the dual of such an object yields a
linear $\sigma(\mathcal{B}^\ast, \mathcal{B}) -
\sigma(\mathcal{A}^\ast, \mathcal{A})$ homeomorphism between the
respective dual spaces. With this observation in place, the proof
of the previous theorem now modifies readily with Proposition
\ref{genbij} being used instead of Theorem \ref{gen}. Note also
that once the existence of $\widetilde{\varphi}$ has been
established in the proof of $(2) \Rightarrow (1)$, it is enough to
verify that$\widetilde{\varphi}(\mathcal{A}) \subset \mathcal{B}$,
since by symmetry we will then also have
$\widetilde{\varphi}^{-1}(\mathcal{B}) \subset \mathcal{A}$
(implying that in this case $\widetilde{\varphi}$ restricts to the
required Jordan $\ast$-isomorphism $\varphi$ from $\mathcal{A}$
onto $\mathcal{B}$).
\end{proof}

\begin{corollary} [Banach-Stone]
A mapping $\varphi$ of $C(K)$ into $C(S)$, with $K$ and $S$ compact Hausdorff
spaces, is a $\ast$-homomorphism if and only if there exists a continuous
transformation $\nu$ of $S$ into $K$ such that
$$\varphi(f) = f \circ \nu \quad \mbox{for each} \quad f \in C(K).$$ In
particular $\varphi$ is a $*$-isomorphism from $C(K)$ onto $C(S)$ if and only
if $\nu$ is a homeomorphism from $S$ onto $K$.
\end{corollary}

\end{document}